\documentclass[11pt,a4paper]{amsart}

\usepackage[a4paper,margin=2.35cm,headheight=15pt]{geometry}
\usepackage[T1]{fontenc}
\usepackage[utf8]{inputenc}
\usepackage{lmodern}
\usepackage{amsmath,amssymb,amsthm,mathtools,mathrsfs}
\usepackage{booktabs,longtable,array,tabularx}
\usepackage{enumerate}
\usepackage{microtype}
\usepackage{hyperref}

\numberwithin{equation}{section}

\newtheorem{theorem}{Theorem}[section]
\newtheorem{lemma}[theorem]{Lemma}
\newtheorem{proposition}[theorem]{Proposition}
\newtheorem{corollary}[theorem]{Corollary}
\theoremstyle{definition}

\theoremstyle{remark}
\newtheorem{remark}[theorem]{Remark}

\title{Sharp Hausdorff--Young inequalities for solvable and semisimple Lie groups}



\author{Jinsong Wu}
\address{Jinsong Wu, Beijing Institute of Mathematical Sciences and Applications, Beijing, 101408, China}
\email{wjs@bimsa.cn}

\author{Fulin Yang}
\address{Fulin Yang, Beijing Institute of Mathematical Sciences and Applications, Beijing, 101408, China}
\address{Yau Mathematical Sciences Center, Tsinghua University, Beijing 100084, China}
\email{fulinyang@bimsa.cn}

\date{}

\begin{document}

\begin{abstract}
In this paper, we obtain the sharp constants for Hausdorff--Young inequalities on simply connected, completely solvable Lie groups, and connected real semisimple Lie groups with finite center.
\end{abstract}

\maketitle


\section{Introduction}\label{sec:main-results}

The sharp Hausdorff–Young inequalities play a fundamental role in harmonic analysis. 
In 1961, Babenko \cite{Babenko1961} proved the sharp Hausdorff–Young inequalities for Euclidean space in special cases, and in 1975, Beckner \cite{Beckner1975} proved them in the general case. 
These inequalities have fruitful consequences in harmonic analysis.

In 1958, Kunze \cite{Kunze1958} proved the Hausdorff–Young inequalities for locally compact unimodular groups. 
In 1980, Terp \cite{Terp2017} proved the Hausdorff–Young inequalities for locally compact groups. 
These Hausdorff–Young inequalities are sharp for compact groups. 
In 2019, Cowling, Martini, Müller, and Parcet \cite{CMMP2019} proved a sharp local central version of the Hausdorff–Young inequalities for compact Lie groups. 
For the Heisenberg group and Weyl transforms, Klein and Russo
\cite{KleinRusso1978} obtained sharp results when the conjugate exponent
is an even integer.  Subsequent estimates for nilpotent Lie groups were
obtained by Baklouti, Smaoui, and Ludwig
\cite{BakloutiSmaouiLudwig2003}.  In the quantum setting, Qiu and Xiong
proved the sharp Hausdorff–Young inequalities for quantum Euclidean
spaces in a recent preprint \cite{QiuXiong2026}.


Let $G$ be a second countable locally compact group with a fixed left Haar measure. 
Denote its left group von Neumann algebra by $\operatorname{VN}(G)$, 
its Fourier algebra by $A(G)=\operatorname{VN}(G)_{*}$, the predual of $\operatorname{VN}(G)$, and its space of bounded left convolution operators by $\operatorname{Cv}^{2}(G)$. 
We use the standard isometric identification $\operatorname{Cv}^{2}(G)=\operatorname{VN}(G)$; see \cite{Eymard1964,Kunze1958}. 
Viewed as spaces of distributions on $G$, $A(G)$ and
$\operatorname{Cv}^{2}(G)$ form a compatible interpolation couple. 
For $1<r<\infty$, set
$$ \mathcal{F}L^{r}(G) = [A(G),\operatorname{Cv}^{2}(G)]_{1-\frac{1}{r}}.$$
At the endpoints, let $\mathcal{F}L^{1}(G)=A(G)$ and
$\mathcal{F}L^{\infty}(G)=\operatorname{Cv}^{2}(G)$. 
The Plancherel identification gives $\mathcal{F}L^{2}(G)=L^{2}(G)$ isometrically. This
interpolation model agrees with the Hilsum--Terp Fourier $L^{r}$-space; see
\cite{Hilsum1981,Terp1982,Terp2017}.
Let $1\leq p\leq2$, let $q=\frac{p}{p-1}$ the H\"{o}lder conjugate of $p$, and define
\begin{align*}
 B_{p}=p^{\frac{1}{2p}}q^{-\frac{1}{2q}},
\end{align*}
with the continuous convention $B_{1}=1$. 
The optimal Hausdorff--Young constant of $G$ is
\begin{align*}
 H_{p}(G) =\sup_{0\neq f\in C_{c}(G)} \frac{\|f\|_{\mathcal{F}L^{q}(G)}}{\|f\|_{L^{p}(G)}}.
\end{align*}
The endpoint estimates and complex interpolation give $H_{p}(G)\leq1$.
If $G$ is unimodular, then
\begin{align*}
 \|f\|_{\mathcal{F}L^{q}(G)}
 =\|\lambda_{G}(f)\|_{L^{q}(\operatorname{VN}(G),\tau_{G})},
\end{align*}
where $\tau_{G}$ is the Plancherel trace, $\lambda_G$ is the left regualar representation. 
If $G$ is nonunimodular, the right-hand side is interpreted in the Hilsum--Terp sense.
In the type-I case, the corresponding Plancherel decomposition of a
nonunimodular group involves the Duflo--Moore operators; see
\cite{DufloMoore1976}.  The Hilsum--Terp formulation used here avoids
the choice of a particular direct-integral realization.

The Babenko--Beckner theorem states that
$$H_{p}(\mathbb{R}^{d})=B_{p}^{d}.$$
While Cowling, Martini, M\"uller, and Parcet proved the universal local lower bound for Lie groups
\begin{align}\label{formula:lower bound} 
H_{p}(G)\geq H_{p}^{\mathrm{loc}}(G)\geq B_{p}^{\dim G},
\end{align}
where $H_{p}^{\mathrm{loc}}(G)$ is the local Hausdorff-Young constant defined in \cite{CMMP2019}.
In this paper,  we shall revise the tricks due to \cite{PotapovSukochevTomskovaZanin2019,QiuXiong2026,CMMP2019,Beckner1975,Izumi1997,Takesaki1973} and prove the sharp Hausdorff-Young inequalities for completely solvable Lie groups and connected real semisimple Lie group with finite center as follows:

\begin{theorem}\label{thm:completely-solvable-sharp}
Let $G$ be a connected, simply connected, completely solvable Lie group.
Then, for every $1\leq p\leq2$,
$$ H_{p}(G)=B_{p}^{\dim G}.$$
\end{theorem}

The examples adapted to Theorem \ref{thm:completely-solvable-sharp} include all connected, simply connected nilpotent Lie groups, the real affine group, real upper-triangular groups with positive diagonal, and the Iwasawa groups $AN$ of connected semisimple Lie groups.  The latter application is developed in the next theorem.
Further nonnilpotent examples include the Damek--Ricci groups, which
are harmonic solvable extensions of $H$-type groups
\cite{DamekRicci1992}.  More generally, a Heintze group whose defining
expanding derivation has only real eigenvalues is completely solvable
and belongs to the class covered by Theorem
\ref{thm:completely-solvable-sharp}; see \cite{Heintze1974}.

\begin{theorem}\label{thm:semisimple-sharp}
Let $G$ be a connected real semisimple Lie group with finite center, and let $K$ be a maximal compact subgroup. Then, for every
$1\leq p\leq2$,
$$H_{p}(G)=B_{p}^{\dim(G/K)}.$$
\end{theorem}
In  Theorem \ref{thm:semisimple-sharp}, the calculation of the sharp Hausdorff--Young constant reduces to the dimension of the associated noncompact Riemannian symmetric space $G/K$.

The global form of a semisimple group requires a small amount of care.  If $\pi:\widetilde{G}\to G$ is a finite central covering and $\widetilde{K}=\pi^{-1}(K)$, then
\begin{align*}
\widetilde{G}/\widetilde{K}\simeq G/K,\qquad \dim(\widetilde{G}/\widetilde{K})=\dim(G/K).
\end{align*}
Thus, all connected global forms with finite center of a fixed real semisimple Lie algebra have the same constant in
Theorem \ref{thm:semisimple-sharp} .  Given a semisimple Lie group, a maximal compact subgroup is therefore understood up to the finite central quotient dictated by the chosen global form.  

Under assumptions of  Theorem \ref{thm:semisimple-sharp}, let $G$ be an almost direct product of noncompact factors $G_{1},\ldots,G_{m}$ with associated maximal compact factors $K_{1},\ldots,K_{m}$. Theorem \ref{thm:semisimple-sharp} and the fact $\dim(G/K)=\sum_{j=1}^{m}\dim(G_{j}/K_{j})$  give the following result.
\begin{corollary}
If a connected finite-center semisimple group $G$ is an almost direct product of noncompact factors $G_{1},\ldots,G_{m}$ with associated maximal compact factors $K_{1},\ldots,K_{m}$, then
\begin{align*}
 H_{p}(G)
 =B_{p}^{\sum_{j=1}^{m}\dim(G_{j}/K_{j})}.
\end{align*}
\end{corollary}

The paper is organized as follows.  Section~\ref{sec:preliminary}
collects the Lie-theoretic, Fourier-analytic, and operator-algebraic
preliminaries, including crossed products, dual weights used throughout the paper.
Section~\ref{sec:one-dimensional-extension} proves the sharp
one-dimensional extension estimate under trace-scaling actions by
combining a weighted crossed-product Kato--Seiler--Simon inequality,
a Hessian inequality in Haagerup spaces, and a flow-and-freezing
argument.  Section~\ref{sec:completely-solvable-proof} iterates this
extension along ideal flags and proves
Theorem~\ref{thm:completely-solvable-sharp}.
Section~\ref{sec:semisimple-proof} proves
Theorem~\ref{thm:semisimple-sharp} through the closed-subgroup
comparison principle, contraction to the Cartan motion group, and the
sharp constant of compact Euclidean extensions.
Finally, Section~\ref{sec:examples} records classical, complex,
compact, and exceptional examples of finite-center semisimple Lie
groups.

\section{Preliminaries}\label{sec:preliminary}

This section collects the structural and operator-algebraic background
used in the proofs.  All locally compact groups are second countable,
all Haar measures are left Haar measures, and all von Neumann algebras
have separable predual.  The statements in this
section are standard and are accompanied by references.

\subsection{Solvable Lie groups}
\label{subsec:preliminary-solvable-Lie-groups}

Let $G$ be a connected Lie group with Lie algebra $\mathfrak{g}$.  The
derived series is
\begin{align*}
 \mathfrak{g}^{(0)}=\mathfrak{g},\qquad
 \mathfrak{g}^{(j+1)}=[\mathfrak{g}^{(j)},\mathfrak{g}^{(j)}],\qquad
 j\geq0.
\end{align*}
The algebra $\mathfrak{g}$, and hence $G$, is called solvable if
$\mathfrak{g}^{(m)}=\{0\}$ for some $m$.  It is called completely
solvable if every eigenvalue of $\operatorname{ad}_{\mathfrak{g}}X$ is
real for every $X\in\mathfrak{g}$.  Complete solvability implies
supersolvability.  We shall use the following standard structure
theorem; see \cite[Corollary~1.30, Corollary~1.126 and
Theorem~1.125]{Knapp2002} and also
\cite{LeptinLudwig1994,Pukanszky1971}.

\begin{lemma}\label{lem:completely-solvable-ideal-flag}
Let $G$ be connected, simply connected, and completely solvable, and
let $n=\dim G$.
\begin{enumerate}[\rm(i)]
\item\label{lem:flag decomposition} There is a complete flag of ideals
\begin{align*}
 \{0\}=\mathfrak{g}_{0}\subset\mathfrak{g}_{1}
 \subset\cdots\subset\mathfrak{g}_{n}=\mathfrak{g},
 \qquad \dim(\mathfrak{g}_{j})=j.
\end{align*}
If $G_{j}$ is the analytic subgroup with Lie algebra
$\mathfrak{g}_{j}$, then $G_{j}$ is closed, connected, simply
connected, and normal in $G$.

\item\label{lem:flag isomorphism} After choosing
$X_{j}\in\mathfrak{g}_{j}\setminus\mathfrak{g}_{j-1}$, there is a
global isomorphism
\begin{align*}
 G_{j}\simeq G_{j-1}\rtimes_{\alpha_{j}}\mathbb{R},\qquad
 \alpha_{j}(t)(x)=\exp(tX_{j})x\exp(-tX_{j}).
\end{align*}
\end{enumerate}
\end{lemma}

For every group in Lemma~\ref{lem:completely-solvable-ideal-flag},
the exponential map
\begin{align*}
 \exp:\mathfrak{g}\longrightarrow G
\end{align*}
is a global diffeomorphism.  More generally, if
$\mathfrak{h}\subseteq\mathfrak{g}$ is an ideal of codimension one,
$H=\exp(\mathfrak{h})$, and
$X_{0}\in\mathfrak{g}\setminus\mathfrak{h}$, then $H$ is closed,
connected, and simply connected, and
\begin{align*}
 H\rtimes_{\alpha}\mathbb{R}&\longrightarrow G,\quad
 (h,t)\longmapsto h\exp(tX_{0})\quad\mbox{with}\quad
 \alpha_{t}(h)=\exp(tX_{0})h\exp(-tX_{0}),
\end{align*}
is a Lie-group isomorphism.  These are standard properties of
simply connected completely solvable groups; see
\cite{LeptinLudwig1994,Pukanszky1971}.

The infinitesimal modular character is
\begin{align*}
 \delta_{\mathfrak{g}}(X)
 =\operatorname{tr}\bigl(\operatorname{ad}_{\mathfrak{g}}X\bigr),
 \qquad X\in\mathfrak{g}.
\end{align*}
Let $\Delta_{G}$ denote the modular function of $G$.  With the
convention used throughout the paper,
\begin{align*}
 \Delta_{G}(\exp X)
 =\exp\bigl(-\delta_{\mathfrak{g}}(X)\bigr).
\end{align*}
Thus $G$ is unimodular if and only if
$\delta_{\mathfrak{g}}=0$.  The additional assertions about ideals and
the modular kernel needed in the solvable case are proved in
Section~\ref{sec:completely-solvable-proof} and are not included among
the preliminaries.

Let $N$ be a unimodular locally compact group and let
$\alpha:\mathbb{R}\rightarrow\operatorname{Aut}(N)$ be a continuous
action.  Its Haar-measure Jacobian is determined by
\begin{align*}
 d(\alpha_{t}(n))=J_{\alpha}(t)\,dn,\qquad
 J_{\alpha}(t)>0.
\end{align*}
The semidirect product $N\rtimes_{\alpha}\mathbb{R}$ has multiplication
and inversion
\begin{align*}
 (n,t)(m,s)=\bigl(n\alpha_{t}(m),t+s\bigr),\qquad
 (n,t)^{-1}=\bigl(\alpha_{-t}(n^{-1}),-t\bigr).
\end{align*}
Its left Haar measure and modular function are
\begin{align*}
 d\mu_{N\rtimes_{\alpha}\mathbb{R}}(n,t)=J_{\alpha}(t)^{-1}\,dn\,dt,\qquad
 \Delta_{N\rtimes_{\alpha}\mathbb{R}}(n,t)=J_{\alpha}(t)^{-1}.
\end{align*}
On
$L^{2}(N\rtimes_{\alpha}\mathbb{R},J_{\alpha}(t)^{-1}\,dn\,dt)$,
\begin{align*}
 \bigl(\lambda_{N\rtimes_{\alpha}\mathbb{R}}(n,t)f\bigr)(m,s)
 =f\bigl(\alpha_{-t}(n^{-1}m),s-t\bigr),
\end{align*}
and the convolution is
\begin{align*}
 (\phi*_{\rtimes_{\alpha}}\psi)(n,t)
 =\int_{\mathbb{R}}\int_{N}
 \phi(m,s)\psi\bigl(\alpha_{-s}(m^{-1}n),t-s\bigr)
 J_{\alpha}(s)^{-1}\,dm\,ds.
\end{align*}
In the trace-scaling situation used below,
$J_{\alpha}(t)=e^{\omega t}$ and hence
\begin{align}\label{eq:nonunimodular-haar}
 d\mu_{G}(n,t)=e^{-\omega t}\,dn\,dt,\qquad
 \Delta_{G}(n,t)=e^{-\omega t},\qquad
 G=N\rtimes_{\alpha}\mathbb{R}.
\end{align}
Consequently, for $1\leq v<\infty$,
$L^{v}(G) =L^{v}\bigl(\mathbb{R},e^{-\omega t}\,dt;L^{v}(N)\bigr)$
isometrically.

\subsection{Semisimple Lie groups}
\label{subsec:preliminary-semisimple-Lie-groups}

Let $G$ be a connected real semisimple Lie group with finite center.
Every connected semisimple Lie group is unimodular, since
$\operatorname{tr}(\operatorname{ad}_{\mathfrak{g}}X)=0$ for every
$X\in\mathfrak{g}$.
Choose a Cartan involution and write
\begin{align*}
 \mathfrak{g}=\mathfrak{k}\oplus\mathfrak{p},\qquad
 [\mathfrak{k},\mathfrak{k}]\subseteq\mathfrak{k},\qquad
 [\mathfrak{k},\mathfrak{p}]\subseteq\mathfrak{p},\qquad
 [\mathfrak{p},\mathfrak{p}]\subseteq\mathfrak{k}.
\end{align*}
The analytic subgroup $K$ with Lie algebra $\mathfrak{k}$ is maximal
compact.  Choose a maximal abelian subspace
$\mathfrak{a}\subseteq\mathfrak{p}$, a positive restricted-root system
$\Sigma^{+}\subseteq\mathfrak{a}^{*}$, and let
\begin{align*}
 \mathfrak{n}=\bigoplus_{\alpha\in\Sigma^{+}}\mathfrak{g}_{\alpha},
 \qquad m_{\alpha}=\dim\mathfrak{g}_{\alpha}.
\end{align*}
Let $A$ and $N$ be the analytic subgroups with Lie algebras
$\mathfrak{a}$ and $\mathfrak{n}$, respectively.
The Iwasawa decompositions \cite{Iwasa1949} are
\begin{align}\label{eq:Iwasawa-decompositions}
 \mathfrak{g}=\mathfrak{k}\oplus\mathfrak{a}\oplus\mathfrak{n},
 \qquad G=KAN,
\end{align}
and multiplication $K\times A\times N\rightarrow G$ is a
diffeomorphism.  The following standard consequence will be used in
the proof of the semisimple theorem; see
\cite[Chapters~VI and IX]{Helgason2001}.

\begin{lemma}\label{lem:AN-completely-solvable}
Under the assumptions of Theorem~\ref{thm:semisimple-sharp} and for
the Iwasawa decomposition~\eqref{eq:Iwasawa-decompositions}, the
subgroup $AN$ is closed, connected, simply connected, and
completely solvable.  Moreover,
\begin{align*}
 \dim(AN)
 =\dim\mathfrak{a}+\sum_{\alpha\in\Sigma^{+}}m_{\alpha}
 =\dim\mathfrak{p}=\dim(G/K).
\end{align*}
\end{lemma}

Choose an $\operatorname{Ad}(K)$-invariant inner product on
$\mathfrak{p}$, let $dX$ be the corresponding Lebesgue measure, and
normalize the Haar measure $dk$ on $K$ by $\int_{K}dk=1$.  The polar
map
\begin{align*}
 \Theta:K\times\mathfrak{p}\longrightarrow G,\qquad
 \Theta(k,X)=k\exp X,
\end{align*}
is a global diffeomorphism.  Normalize the Haar measure $dg$ on $G$
so that there is a unique positive smooth
$\operatorname{Ad}(K)$-invariant function $J$ satisfying
\begin{align}\label{eq:Cartan-Jacobian}
 \Theta^{*}(dg)=J(X)\,dk\,dX,\qquad J(0)=1.
\end{align}
These facts follow from the Cartan and polar decompositions; see
\cite[Chapters~VI and IX]{Helgason2001}.

More generally, let $V$ be a finite-dimensional Euclidean space and
let $\rho:K\rightarrow O(V)$ be an orthogonal representation of a
compact group $K$.  Write $k\mathbin{\cdot}x=\rho(k)x$.  The standard
semidirect product $K\ltimes V$ is $K\times V$ with
\begin{align*}
 (k,x)(l,y)=(kl,x+k\mathbin{\cdot}y),\qquad
 (k,x)^{-1}=(k^{-1},-k^{-1}\mathbin{\cdot}x).
\end{align*}
The standard Cartan motion group is
$K\ltimes_{\operatorname{Ad}}\mathfrak{p}$.  A function
$g:V\rightarrow\mathbb{C}$ is $K$-invariant if
\begin{align*}
 g(k\mathbin{\cdot}x)=g(x),\qquad k\in K,\qquad x\in V.
\end{align*}

We also recall the quotient disintegration used in the closed-subgroup
argument.  For a closed subgroup $H$ of a locally compact group $G$,
write $H_{G}=\{Hg:g\in G\}$ for the left coset space.

\begin{lemma}\label{lem:left-coset-disintegration}
Let $G$ be a second countable locally compact group and let $H$ be a
closed subgroup.  After fixing left Haar measures on $G$ and $H$, the
following statements hold.
\begin{enumerate}[\rm(i)]
\item\label{coset representation} There are a Borel section
$s:H_{G}\rightarrow G$ and a $\sigma$-finite Borel measure $\nu$ on
$H_{G}$ such that
\begin{align*}
 \int_{G}F(g)\,dg
 =\int_{H_{G}}\int_{H}F(hs(x))\,dh\,d\nu(x)
\end{align*}
for every nonnegative Borel function $F$ on $G$.

\item\label{unitary adjoint represation} The operator
\begin{align*}
 \mathcal{U}:L^{2}(G)\longrightarrow L^{2}(H_{G};L^{2}(H)),\qquad
 (\mathcal{U}\xi)(x,h)=\xi(hs(x)),
\end{align*}
is unitary and satisfies
\begin{align*}
 \mathcal{U}\lambda_{G}(h_{0})\mathcal{U}^{-1}
 =I_{L^{2}(H_{G})}\otimes\lambda_{H}(h_{0}),\qquad h_{0}\in H.
\end{align*}
\end{enumerate}
\end{lemma}

Lemma~\ref{lem:left-coset-disintegration} is the standard measurable
disintegration over the left coset space.  See
\cite{Folland2016}, Sections~2.6--2.7.  No invariance of $\nu$ under
$G$ is asserted or needed.

\subsection{Fourier transform and group von Neumann algebras}
\label{subsec:preliminary-group-Fourier-transform}

For scalar functions on $\mathbb{R}$, we use the normal Fourier transform
\begin{align}\label{eq:scalar-Fourier-transform}
 \widehat{b}(\xi)=\int_{\mathbb{R}}b(t)e^{-2\pi i t\xi}\,dt.
\end{align}

Let $G$ be a second countable locally compact group with left Haar
measure $dx$ and modular function $\Delta_{G}$.  Its left regular
representation on $L^{2}(G)$ is
\begin{align*}
 (\lambda_{G}(x)\xi)(y)=\xi(x^{-1}y),\qquad x,y\in G.
\end{align*}
For $f\in L^{1}(G)$, let
\begin{align*}
 \lambda_{G}(f)=\int_{G}f(x)\lambda_{G}(x)\,dx.
\end{align*}
Then $\|\lambda_{G}(f)\|_{\operatorname{VN}(G)} \leq\|f\|_{L^{1}(G)},$
where
\begin{align*}
 \operatorname{VN}(G)
 =\{\lambda_{G}(x):x\in G\}^{\prime\prime}.
\end{align*}
The predual of $\operatorname{VN}(G)$ is the Fourier algebra $A(G)$.
We identify $\operatorname{VN}(G)$ with the bounded left convolution
operators $\operatorname{Cv}^{2}(G)$; see
\cite{Eymard1964,Kunze1958}.

The left Hilbert algebra of $G$ determines a normal faithful semifinite
Plancherel weight $\varphi_{G}$ on $\operatorname{VN}(G)$, characterized
on its natural convolution core by
\begin{align*}
 \varphi_{G}\bigl(\lambda_{G}(f)^{*}\lambda_{G}(f)\bigr)
 =\|f\|_{L^{2}(G)}^{2}.
\end{align*}
It is tracial exactly when $G$ is unimodular.  The Fourier
$L^{r}$-spaces are
\begin{align*}
 \mathcal{F}L^{r}(G)
 =[A(G),\operatorname{Cv}^{2}(G)]_{1-\frac{1}{r}},\qquad
 1\leq r\leq\infty,
\end{align*}
and, by reiteration,
\begin{align}\label{FL interpolation 2-infty}
 \mathcal{F}L^{r}(G)
 =[\operatorname{Cv}^{2}(G),L^{2}(G)]_{\frac{2}{r}},\qquad
 2\leq r\leq\infty.
\end{align}
These descriptions agree with the Hilsum--Terp construction; see
\cite{Hilsum1981,Terp1982,Terp2017,Izumi1997}.

Let $\mathcal{H}_{\varphi_{G}}$ be the GNS Hilbert space of
$\varphi_{G}$.  There are canonical complex-linear isometries
\begin{align*}
 \mathcal{F}L^{r}(G)
 \cong
 [\operatorname{VN}(G),\mathcal{H}_{\varphi_{G}}]_{\frac{2}{r}}
 =L^{r}(\operatorname{VN}(G))_{\mathrm{left}},\qquad
 2\leq r\leq\infty.
\end{align*}
On the Fourier--Tomita core, let $j_{r}^{\varphi_{G}}$ denote the
concrete interpolation map.  Thus
\begin{align*}
 \|f\|_{\mathcal{F}L^{r}(G)}
 =\|j_{r}^{\varphi_{G}}(\lambda_{G}(f))\|_{L^{r}(\operatorname{VN}(G))_{\mathrm{left}}}.
\end{align*}
If $G$ is unimodular, then $\varphi_{G}=\tau_{G}$ and this is the usual
tracial space $L^{r}(\operatorname{VN}(G),\tau_{G})$.

For $1\leq p\leq2$ and $q=\frac{p}{p-1}$, the interpolation extension
of $f\mapsto\lambda_{G}(f)$ is denoted by
\begin{align*}
 \mathcal{F}_{G,p}:L^{p}(G)\longrightarrow\mathcal{F}L^{q}(G).
\end{align*}
On the Fourier--Tomita core,
\begin{align*}
 \mathcal{F}_{G,p}f\longleftrightarrow j_{q}^{\varphi_{G}}(\lambda_{G}(f)).
\end{align*}

The standard Tomita-core facts will be used in both the solvable and
semisimple arguments.

\begin{lemma}\label{lem:semisimple-Tomita-core}
Let $G$ be a second countable Lie group.  The convolution algebra
$C_{c}^{\infty}(G)$, with involution
\begin{align*}
 f^{\sharp}(x)
 =\Delta_{G}(x)^{-1}\overline{f(x^{-1})},
\end{align*}
is an entire analytic function-side Tomita core for the Plancherel
weight.  Let $\sigma^{\varphi_{G}}$ denote its modular automorphism
group.  On this core,
\begin{align*}
 \sigma_{z}^{\varphi_{G}}(\lambda_{G}(f))
 =\lambda_{G}(\Delta_{G}^{iz}f),\qquad z\in\mathbb{C}.
\end{align*}
Moreover, $C_{c}^{\infty}(G)$ is dense in $L^{v}(G)$ for
$1\leq v<\infty$, and it is dense in $\mathcal{F}L^{r}(G)$ for
$2\leq r<\infty$.  In addition,
\begin{align*}
 [\operatorname{Cv}^{2}(G),L^{2}(G)]_{\theta}
 =\mathcal{F}L^{\frac{2}{\theta}}(G),\qquad 0<\theta<1,
\end{align*}
isometrically.
\end{lemma}

See \cite[Chapter~VIII]{Takesaki2003} and
\cite{Hilsum1981,Terp1982,Terp2017,Izumi1997} for
Lemma~\ref{lem:semisimple-Tomita-core}.

\subsection{Crossed products}
\label{subsec:preliminary-crossed-products}

Let $\mathcal{M}$ be a von Neumann algebra acting standardly on a
Hilbert space $\mathcal{H}$, and let
$\beta:\mathbb{R}\rightarrow\operatorname{Aut}(\mathcal{M})$ be a
pointwise weak-* continuous action.  On
$L^{2}(\mathbb{R};\mathcal{H})$, let
\begin{align*}
 (\pi_{\beta}(x)\xi)(s)=\beta_{-s}(x)\xi(s),\qquad
 (\lambda_{t}\xi)(s)=\xi(s-t).
\end{align*}
The covariance relation is 
\begin{align*}
 \lambda_{t}\pi_{\beta}(x)\lambda_{t}^{*}
 =\pi_{\beta}(\beta_{t}(x)).
\end{align*}
The crossed product is von Neumann algebra
\begin{align*}
 \mathcal{M}\rtimes_{\beta}\mathbb{R}
 =\bigl(\pi_{\beta}(\mathcal{M})
 \cup\{\lambda_{t}:t\in\mathbb{R}\}\bigr)^{\prime\prime}.
\end{align*}
For the crossed-product construction, its duality theory, and the
associated dual weights, see \cite{Takesaki1973}.
For a compactly supported coefficient function
$F:\mathbb{R}\rightarrow\mathcal{M}$, let
\begin{align}\label{eq:definition-TF}
 T_{F}=\int_{\mathbb{R}}\pi_{\beta}(F(t))\lambda_{t}\,dt.
\end{align}
The twisted convolution and involution are
\begin{align*}
 (F_{1}*_{\beta}F_{2})(t)=\int_{\mathbb{R}}F_{1}(s)\beta_{s}(F_{2}(t-s))\,ds,\qquad
 F^{\sharp}(t)=\beta_{t}(F(-t)^{*}).
\end{align*}
Thus $T_{F_{1}}T_{F_{2}}=T_{F_{1}*_{\beta}F_{2}}$ and
$T_{F}^{*}=T_{F^{\sharp}}$ on the compactly supported left Hilbert
algebra.

\subsection{Dual weights}
\label{subsec:preliminary-dual-weights}

Suppose that $\mathcal{M}$ carries a normal faithful semifinite trace
$\tau$ and, for some $\omega\in\mathbb{R}$,
\begin{align*}
 \tau\circ\beta_{t}=e^{-\omega t}\tau,\qquad t\in\mathbb{R}.
\end{align*}
In all applications below, $\omega\geq0$. 
On $\mathcal{M}\rtimes_{\beta}\mathbb{R}$, the dual action is determined by
\begin{align*}
 \widehat{\beta}_{s}(\pi_{\beta}(x))=\pi_{\beta}(x),\qquad
 \widehat{\beta}_{s}(\lambda_{t})=e^{ist}\lambda_{t}.
\end{align*}
With this convention, the operator-valued weight $\mathcal{T}_{\widehat{\beta}}:
 (\mathcal{M}\rtimes_{\beta}\mathbb{R})_{+}\longrightarrow\widehat{\pi_{\beta}(\mathcal{M})}_{+}$ is given by
\begin{align*}
 \mathcal{T}_{\widehat{\beta}}(X)
 =\int_{\mathbb{R}}\widehat{\beta}_{s}(X)\,\frac{ds}{2\pi},
 \qquad X\in(\mathcal{M}\rtimes_{\beta}\mathbb{R})_{+},
\end{align*}
where the integral is taken in the extended positive cone.  The dual
weight associated with $\tau$ is, by definition,
\begin{align*}
 \Phi
 =\tau\circ\pi_{\beta}^{-1}\circ
 \mathcal{T}_{\widehat{\beta}}.
\end{align*}
Here $\pi_{\beta}^{-1}$ denotes its natural extension to the extended
positive cone of $\pi_{\beta}(\mathcal{M})$.
It is normal, faithful, semifinite, and invariant under the dual
action.  Its square-integrable left ideal and GNS map are
$$\mathfrak{n}_{\Phi}=\{Y\in\mathcal{M}\rtimes_{\beta}\mathbb{R}:\Phi(Y^{*}Y)<\infty\}$$
and
$$ \Lambda_{\Phi}:\mathfrak{n}_{\Phi}\longrightarrow\mathcal{H}_{\Phi}
 \quad\mbox{with}\quad
 \|\Lambda_{\Phi}(Y)\|_{\mathcal{H}_{\Phi}}^{2}=\Phi(Y^{*}Y).$$
We use the standard notation
$\mathfrak{n}_{\Phi}^{*}=\{Y\in\mathcal{M}\rtimes_{\beta}\mathbb{R}:Y^{*}\in
\mathfrak{n}_{\Phi}\}$.
Let $\mathfrak{C}_{\beta}$ be the algebra generated, under
$*_{\beta}$ and $\sharp$, by coefficient functions
\begin{align*}
 F(t)=b(t)x,\qquad
 b\in C_{c}^{\infty}(\mathbb{R}),\qquad
 x\in\mathcal{M}\cap L^{2}(\mathcal{M},\tau).
\end{align*}
The corresponding coefficient Tomita core is
$\{T_{F}:F\in\mathfrak{C}_{\beta}\}
 \subseteq\mathfrak{n}_{\Phi}\cap\mathfrak{n}_{\Phi}^{*}.$ 
The explicit GNS inner-product formula on this core is proved in
Lemma~\ref{lem:trace-scaling-dual-weight}; it is not part of the
definition of $\Phi$.

Since $\tau$ is tracial, the modular group of the dual weight satisfies
\begin{align*}
 \sigma_{s}^{\Phi}(\pi_{\beta}(x))=\pi_{\beta}(x),\qquad
 \sigma_{s}^{\Phi}(\lambda_{t})=e^{-i\omega ts}\lambda_{t}.
\end{align*}
Let $D_{\Phi}$ denote the canonical positive self-adjoint density of
$\Phi$ in the Haagerup core, so that
$D_{\Phi}^{is}YD_{\Phi}^{-is}=\sigma_{s}^{\Phi}(Y)$.  For $r>0$, on
the analytic core,
\begin{align}\label{eq:dual-density-commutation}
 D_{\Phi}^{1/r}\lambda_{t}
 =e^{-\frac{\omega t}{r}}\lambda_{t}D_{\Phi}^{1/r}.
\end{align}

In the left Izumi--Hilsum--Terp compatible realization, an element
$Y\in\mathfrak{n}_{\Phi}$ is represented by $Y$ at the
$\mathcal{M}\rtimes_{\beta}\mathbb{R}$-endpoint and by $\Lambda_{\Phi}(Y)$ at the Hilbert-space
endpoint.  Under these canonical compatible embeddings, the common
core is identified with $\mathfrak{n}_{\Phi}$.  For
$2\leq r<\infty$,
\begin{align}\label{eq:izumi-left-interpolation}
 [\mathcal{M}\rtimes_{\beta}\mathbb{R},\mathcal{H}_{\Phi}]_{\frac{2}{r}}
 =L^{r}(\mathcal{M}\rtimes_{\beta}\mathbb{R})_{\mathrm{left}}
\end{align}
isometrically.  On the coefficient Tomita core,
$j_{r}^{\Phi}(Y)=\overline{YD_{\Phi}^{1/r}}.$ 
We let $j_{\infty}^{\Phi}(Y)=Y$; at $r=2$,
$j_{2}^{\Phi}(Y)$ corresponds to $\Lambda_{\Phi}(Y)$.  Similarly,
\begin{align*}
 [\mathcal{M},L^{2}(\mathcal{M},\tau)]_{\frac{2}{r}}
 =L^{r}(\mathcal{M},\tau),
\end{align*}
and $j_{r}^{\tau}(x)$ denotes the tracial interpolation image of
$x\in\mathcal{M}\cap L^{2}(\mathcal{M},\tau)$.

These standard facts are established in
\cite{Hilsum1981,Terp1982,Izumi1997,Caspers2013,Takesaki1973,Takesaki2003}; see also
\cite{Kosaki1984,PisierXu2003} for noncommutative interpolation.  If $\omega=0$,
then $\Phi$ is the canonical trace.

\section{One-dimensional extensions under trace-scaling actions}\label{sec:one-dimensional-extension}

This section gives a single proof for the trace-preserving and nonunimodular one-dimensional extensions.  The
essential new ingredient is a crossed-product Kato--Seiler--Simon estimate in the left Izumi--Hilsum--Terp
$L^{q}$-space.
Related Hausdorff--Young estimates for group extensions, including
nonunimodular reductions to the modular kernel under type-I hypotheses,
were obtained by F\"uhr \cite{Fuhr2006}.

\begin{theorem}\label{thm:one-dimensional-extension}
Let $N$ be a simply connected, unimodular, completely solvable Lie group, and let $\alpha:\mathbb{R}\rightarrow\operatorname{Aut}(N)$ be a smooth action such that $N\rtimes_{\alpha}\mathbb{R}$ is completely solvable with $|\det((d\alpha_{t})_{e})|=e^{\omega t}$ for some $\omega\geq0$. For every $1\leq p\leq2$, we have
$$ H_{p}(N\rtimes_{\alpha}\mathbb{R})\leq B_{p}H_{p}(N).$$
\end{theorem}
When $\omega=0$, the formula gives the unimodular product measure.

\subsection{The crossed-product Kato--Seiler--Simon estimate}

We use the coefficient notation $T_{F}$ and the coefficient Tomita core
introduced in Sections~\ref{subsec:preliminary-crossed-products}
and~\ref{subsec:preliminary-dual-weights}, respectively.

\begin{lemma}\label{lem:trace-scaling-dual-weight}
Let $\mathcal{M}$ be a von Neumann algebra with normal faithful semifinite trace $\tau$, and let $\Phi$ be the dual weight associated with $\tau$ on the crossed product $\mathcal{M}\rtimes_{\beta}\mathbb{R}$. Suppose $\tau\circ\beta_{t}=e^{-\omega t}\tau$ for $t\in\mathbb{R}$ and some $\omega\geq0$. If $F_{1},F_{2}\in C_{c}\bigl(\mathbb{R};\mathcal{M}\cap L^{2}(\mathcal{M},\tau)\bigr)$ belong to the coefficient Tomita core, then
$$ \Phi(T_{F_{1}}^{*}T_{F_{2}})=\int_{\mathbb{R}} e^{\omega t} \tau(F_{1}(t)^{*}F_{2}(t))\,dt.$$
\end{lemma}
\begin{proof}
Setting $F^{\sharp}(t)=\beta_{t}(F(-t)^{*})$, we have
\begin{align*}
(F_{1}^{\sharp}*_{\beta}F_{2})(0)
=\int_{\mathbb{R}}\beta_{s}\bigl(F_{1}(-s)^{*}F_{2}(-s)\bigr)\,ds.
\end{align*}
The dual-weight formula evaluates a coefficient at zero.  Hence, with
$\Phi(T_{F_{1}}^{*}T_{F_{2}})$ understood as the sesquilinear extension
of $\Phi$ on its square-integrable ideal, the trace-scaling identity gives
\begin{align*}
 \Phi(T_{F_{1}}^{*}T_{F_{2}})
 =\int_{\mathbb{R}}e^{-\omega s} \tau(F_{1}(-s)^{*}F_{2}(-s))\,ds
 =\int_{\mathbb{R}}e^{\omega t}\tau(F_{1}(t)^{*}F_{2}(t))\,dt.
\end{align*}
The assertion follows.
\end{proof}

By Stone's theorem, let $A$ be the self-adjoint generator of
$\{\lambda_{t}\}_{t\in\mathbb{R}}$, normalized by
$\lambda_{t}=e^{-2\pi itA}.$
The operator $A$ is affiliated with
$\{\lambda_{t}:t\in\mathbb{R}\}^{\prime\prime}\subseteq\mathcal{M}\rtimes_{\beta}\mathbb{R}$.
For $a\in L^{\infty}(\mathbb{R})$, spectral calculus gives
\begin{align*}
 a(A)\in\mathcal{M}\rtimes_{\beta}\mathbb{R},\qquad
 \|a(A)\|_{\mathcal{M}\rtimes_{\beta}\mathbb{R}}=\|a\|_{L^{\infty}(\mathbb{R})}.
\end{align*}
If $b\in C_{c}^{\infty}(\mathbb{R})$, then
$\widehat{b}(A)=\int_{\mathbb{R}}b(t)\lambda_{t}\,dt.$
Combining this identity with
\eqref{eq:dual-density-commutation}, one obtains on the common analytic
core, for $x\in\mathcal{M}\cap L^{2}(\mathcal{M},\tau)$ and
$b\in C_{c}^{\infty}(\mathbb{R})$,
\begin{align*}
 \overline{D_{\Phi}^{1/r}\pi_{\beta}(x)\widehat{b}(A)}
 =\overline{\left(\int_{\mathbb{R}}e^{-\frac{\omega t}{r}}b(t)
 \pi_{\beta}(x)\lambda_{t}\,dt\right)D_{\Phi}^{1/r}}
 =j_{r}^{\Phi}\left(
 \int_{\mathbb{R}}e^{-\frac{\omega t}{r}}b(t)
 \pi_{\beta}(x)\lambda_{t}\,dt\right).
\end{align*}

\begin{lemma}\label{lem:weighted-kss-common-analytic-family}
Let $x\in\mathcal{M}\cap L^{2}(\mathcal{M},\tau)$ and $b\in C_{c}^{\infty}(\mathbb{R})$. For $0\leq\operatorname{Re}(z)\leq1$, define
$$ \mathcal{K}_{z}(x,\widehat{b}) =\int_{\mathbb{R}}e^{-\frac{\omega t z}{2}}b(t)\pi_{\beta}(x)\lambda_{t}\,dt,$$
where the Fourier transform is given in \eqref{eq:scalar-Fourier-transform}. Then
\begin{align*}
 \mathcal{K}_{z}(x,\widehat{b})\in(\mathcal{M}\rtimes_{\beta}\mathbb{R})\cap\mathcal{H}_{\Phi}
\end{align*}
throughout the closed strip, and the mapping
$z\mapsto\mathcal{K}_{z}(x,\widehat{b})$ is analytic in the intersection norm.
Moreover, uniformly for $y\in\mathbb{R}$,
\begin{align*}
 \|\mathcal{K}_{iy}(x,\widehat{b})\|_{\mathcal{M}\rtimes_{\beta}\mathbb{R}}
 \leq\|x\|_{\mathcal{M}}\|\widehat{b}\|_{L^{\infty}(\mathbb{R})},\quad
 \|\Lambda_{\Phi}(\mathcal{K}_{1+iy}(x,\widehat{b}))\|_{\mathcal{H}_{\Phi}}
 =\|x\|_{L^{2}(\mathcal{M},\tau)}\|\widehat{b}\|_{L^{2}(\mathbb{R})}.
\end{align*}
\end{lemma}
\begin{proof}
Since $b$ has compact support, differentiation with respect to $z$ only
introduces the factors
\begin{align*}
 \left(-\frac{\omega t}{2}\right)^{m}e^{-\frac{\omega t z}{2}}b(t),
 \qquad m\geq0.
\end{align*}
They are uniformly integrable on compact subsets of the strip.  The
integrated representation is bounded by the $L^{1}$-norm of its
coefficient.  Taking the difference quotients of
$e^{-\omega t z/2}b(t)$ in $L^{1}(\mathbb{R})$ therefore proves operator-norm
analyticity; no operator-norm continuity of $t\mapsto\lambda_{t}$ is
used.  By Lemma~\ref{lem:trace-scaling-dual-weight},
\begin{align*}
 \|\Lambda_{\Phi}(\mathcal{K}_{z}(x,\widehat{b}))\|_{\mathcal{H}_{\Phi}}^{2}
 = \|x\|_{L^{2}(\mathcal{M},\tau)}^{2}\int_{\mathbb{R}}e^{\omega t(1-\operatorname{Re}(z))}|b(t)|^{2}\,dt,
\end{align*}
and the same dominated-convergence argument applied to difference
quotients proves analyticity with values in $\mathcal{H}_{\Phi}$.  Thus
the family is analytic in the intersection norm.

For $z=iy$, denote $b_{y}(t)=e^{-\frac{i\omega y t}{2}}b(t)$ for $t,y\in\mathbb{R}$. We have
\begin{align*}
\|\mathcal{K}_{iy}(x,\widehat{b})\|_{\mathcal{M}\rtimes_{\beta}\mathbb{R}}
=\|\pi_{\beta}(x)\widehat{b_{y}}(A)\|_{\mathcal{M}\rtimes_{\beta}\mathbb{R}}
\leq\|\pi_{\beta}(x)\|_{\mathcal{M}\rtimes_{\beta}\mathbb{R}}\|\widehat{b_{y}}(A)\|_{\mathcal{M}\rtimes_{\beta}\mathbb{R}}
\leq\|x\|_{\mathcal{M}}\|\widehat{b}\|_{L^{\infty}(\mathbb{R})}.
\end{align*}
Here $\widehat{b_{y}}$ is a translate of $\widehat{b}$, so their
$L^{\infty}$-norms agree.  On the other boundary, Lemma
\ref{lem:trace-scaling-dual-weight} and the Plancherel identity give
\begin{align*}
 \|\Lambda_{\Phi}(\mathcal{K}_{1+iy}(x,\widehat{b}))\|_{\mathcal{H}_{\Phi}}^{2}
 = \int_{\mathbb{R}} e^{\omega t}e^{-\omega t}|b(t)|^{2}\tau(x^{*}x)\,dt=\|x\|_{2}^{2}\|\widehat{b}\|_{2}^{2}.
\end{align*}
The proof is complete.
\end{proof}

\begin{theorem}\label{thm:weighted-crossed-product-kss}
Let $2\leq r<\infty$.  The bilinear assignment initially defined by
\begin{align*}
 \left(j_{r}^{\tau}(x),\widehat{b}\right)
 \longmapsto
 j_{r}^{\Phi}\left(
 \int_{\mathbb{R}}e^{-\frac{\omega t}{r}}b(t)\pi_{\beta}(x)\lambda_{t}\,dt
 \right),
\end{align*}
where $x\in\mathcal{M}\cap L^{2}(\mathcal{M},\tau)$ and
$b\in C_{c}^{\infty}(\mathbb{R})$, extends uniquely to a contraction
$L^{r}(\mathcal{M},\tau)\times L^{r}(\mathbb{R})
 \longrightarrow
 L^{r}(\mathcal{M}\rtimes_{\beta}\mathbb{R})_{\mathrm{left}}.$
In particular, on this core,
\begin{align*}
 \left\|j_{r}^{\Phi}\left(
 \int_{\mathbb{R}}e^{-\frac{\omega t}{r}}b(t)\pi_{\beta}(x)\lambda_{t}\,dt
 \right)\right\|_{L^{r}(\mathcal{M}\rtimes_{\beta}\mathbb{R})_{\mathrm{left}}}
 \leq \|j_{r}^{\tau}(x)\|_{L^{r}(\mathcal{M},\tau)}
 \|\widehat{b}\|_{L^{r}(\mathbb{R})}.
\end{align*}
For general $L^{r}$ inputs, the value is understood as this continuous
extension.  For $r=\infty$,
$x\in\mathcal{M}$, and $\widehat{b}\in L^{\infty}(\mathbb{R})$, one has
\begin{align*}
 \|\pi_{\beta}(x)\widehat{b}(A)\|_{\mathcal{M}\rtimes_{\beta}\mathbb{R}}
 \leq \|x\|_{\mathcal{M}}\|\widehat{b}\|_{L^{\infty}(\mathbb{R})}.
\end{align*}
Here $\widehat{b}(A)$ is defined by spectral calculus; the time-integral
formula is used literally only on the coefficient core.
\end{theorem}
\begin{proof}
Let $\theta=\frac{2}{r}$.  Lemma
\ref{lem:weighted-kss-common-analytic-family} and the bilinear analytic
interpolation theorem \cite[Section~4.4]{BerghLofstrom1976} apply to the
compatible couples below; see also \cite{FournierRusso1977} for
abstract interpolation of operator-valued kernels:
\begin{align*}
 (\mathcal{M},L^{2}(\mathcal{M},\tau)),\qquad
 (L^{\infty}(\mathbb{R}),L^{2}(\mathbb{R})),\qquad
 (\mathcal{M}\rtimes_{\beta}\mathbb{R},\mathcal{H}_{\Phi}).
\end{align*}
More explicitly, on dense admissible analytic test functions in the two
source couples, pair the analytic function obtained by applying
$\mathcal{K}_{z}$ with an admissible analytic test function from the
dual target couple and apply the three-lines theorem.  This is the
analytic-family version of the bilinear interpolation argument used
here.
The interpolation identities  are
\begin{align*}
 [\mathcal{M},L^{2}(\mathcal{M},\tau)]_{\theta}
 =L^{r}(\mathcal{M},\tau),\quad
 [L^{\infty}(\mathbb{R}),L^{2}(\mathbb{R})]_{\theta}
 =L^{r}(\mathbb{R}),\quad
 [\mathcal{M}\rtimes_{\beta}\mathbb{R},\mathcal{H}_{\Phi}]_{\theta}
 =L^{r}(\mathcal{M}\rtimes_{\beta}\mathbb{R})_{\mathrm{left}}.
\end{align*}
At $z=\theta$, the coefficient in the analytic family is
$e^{-\omega t z/2}=e^{-\omega t/r}$, so interpolation gives the asserted core
estimate.  The space $\mathcal{M}\cap L^{2}(\mathcal{M},\tau)$ is dense
in $L^{r}(\mathcal{M},\tau)$, and the Fourier transforms of
$C_{c}^{\infty}(\mathbb{R})$ are dense in $L^{r}(\mathbb{R})$; see
\cite{BerghLofstrom1976,PisierXu2003}.
This proves the unique extension.  The $r=\infty$ estimate follows
directly from spectral calculus.
\end{proof}

\begin{remark}
The proof must use the genuine compatible couple $(\mathcal{M}\rtimes_{\beta}\mathbb{R},\mathcal{H}_{\Phi})$.  A nonzero homogeneous Haagerup $L^{r}(\mathcal{M}\rtimes_{\beta}\mathbb{R})$ element generally has infinite ordinary $L^{r}$--norm with respect to the semifinite trace of the continuous core.  Therefore an unnormalized continuous-core cutoff followed by ordinary tracial interpolation does not give the extension in Theorem \ref{thm:weighted-crossed-product-kss} for $2\leq r<\infty$.
\end{remark}

When $\omega=0$, all exponential weights above equal one and $\Phi$ is the canonical trace. Consequently, Theorem \ref{thm:weighted-crossed-product-kss} reduces to the trace-preserving crossed-product Kato--Seiler--Simon inequality with constant one.

\begin{lemma}\label{lem:weighted-crossed-product-assembly}
Let $1\leq p\leq2$ and $q=\frac{p}{p-1}$, with $q=\infty$ when
$p=1$.  For
$X\in C_{c}\bigl(\mathbb{R};\mathcal{M}\cap L^{2}(\mathcal{M},\tau)\bigr),$
define on the common coefficient core
\begin{align*}
 \mathcal{A}_{\omega}\bigl(t\longmapsto j_{q}^{\tau}(X(t))\bigr)
 =j_{q}^{\Phi}\left(
 \int_{\mathbb{R}}e^{-\omega t}\pi_{\beta}(X(t))\lambda_{t}\,dt
 \right).
\end{align*}
Then $\mathcal{A}_{\omega}$ extends uniquely to a contraction
$L^{p}\bigl(\mathbb{R},e^{-\omega t}\,dt;L^{q}(\mathcal{M},\tau)\bigr)
 \longrightarrow
 L^{q}(\mathcal{M}\rtimes_{\beta}\mathbb{R})_{\mathrm{left}}.$
In particular,
\begin{align*}
 \|\mathcal{A}_{\omega}(X)\|_{L^{q}(\mathcal{M}\rtimes_{\beta}\mathbb{R})_{\mathrm{left}}}
 \leq\left(\int_{\mathbb{R}}e^{-\omega t}\|X(t)\|_{L^{q}(\mathcal{M},\tau)}^{p}\,dt\right)^{\frac{1}{p}}.
\end{align*}
For a general $L^{q}$-valued function, $\mathcal{A}_{\omega}(X)$ denotes
the continuous extension.
\end{lemma}
\begin{proof}
At the $p=1$ endpoint, define the map on the whole space
$L^{1}(\mathbb{R},e^{-\omega t}\,dt;\mathcal{M})$ by the ultraweak integral
\begin{align*}
 \mathcal{A}_{\omega}(X)
 =\int_{\mathbb{R}}^{\mathrm{uw}}e^{-\omega t}
 \pi_{\beta}(X(t))\lambda_{t}\,dt.
\end{align*}
Strong measurability of $X$ implies ultraweak measurability of the
integrand, and evaluation against normal functionals gives
\begin{align*}
 \left\|\int_{\mathbb{R}}^{\mathrm{uw}}e^{-\omega t}
 \pi_{\beta}(X(t))\lambda_{t}\,dt\right\|_{\mathcal{M}\rtimes_{\beta}\mathbb{R}}
 \leq\int_{\mathbb{R}}e^{-\omega t}\|X(t)\|_{\mathcal{M}}\,dt.
\end{align*}
Thus no operator-norm density of
$\mathcal{M}\cap L^{2}(\mathcal{M},\tau)$ in $\mathcal{M}$ is used.
At the $p=2$ endpoint, Lemma
\ref{lem:trace-scaling-dual-weight} gives, on the common coefficient core,
\begin{align*}
 \left\|\Lambda_{\Phi}\left(\int_{\mathbb{R}}e^{-\omega t}
 \pi_{\beta}(X(t))\lambda_{t}\,dt\right)\right\|_{\mathcal{H}_{\Phi}}^{2}
 =\int_{\mathbb{R}}e^{-\omega t}\|X(t)\|_{L^{2}(\mathcal{M},\tau)}^{2}\,dt.
\end{align*}
This map therefore extends isometrically to the whole $L^{2}$ endpoint.
The two endpoint maps agree on the coefficient Tomita core as elements
of the compatible couple
$\bigl(\mathcal{M}\rtimes_{\beta}\mathbb{R},\mathcal{H}_{\Phi}\bigr)$.

Set $\theta=\frac{2}{q}$.  Letting $\frac{1}{p}=1-\frac{\theta}{2}$, vector-valued complex interpolation \cite{BerghLofstrom1976,PisierXu2003} gives
\begin{align*}
 \left[
 L^{1}\bigl(\mathbb{R},e^{-\omega t}\,dt;\mathcal{M}\bigr),
 L^{2}\bigl(\mathbb{R},e^{-\omega t}\,dt;
 L^{2}(\mathcal{M},\tau)\bigr)
 \right]_{\theta}
 &=
 L^{p}\left(
 \mathbb{R},e^{-\omega t}\,dt;
 [\mathcal{M},L^{2}(\mathcal{M},\tau)]_{\theta}
 \right)\\
 &=
 L^{p}\bigl(\mathbb{R},e^{-\omega t}\,dt;
 L^{q}(\mathcal{M},\tau)\bigr).
\end{align*}
Thus the interpolated source is the space in the statement, while
\eqref{eq:izumi-left-interpolation} identifies the target space with
$L^{q}(\mathcal{M}\rtimes_{\beta}\mathbb{R})_{\mathrm{left}}$.
Interpolation proves the contraction.  On the common core, the
interpolated source element is $t\mapsto j_{q}^{\tau}(X(t))$, and the
target element is precisely the $j_{q}^{\Phi}$-expression in the
statement.  This also proves uniqueness by density for $1<p\leq2$;
the $p=1$ endpoint was defined directly on its full source space.
\end{proof}

\subsection{The Hessian inequality in Haagerup spaces}
The flow uses a second-order inequality for the norm on a possibly
nonsemifinite von Neumann algebra.  We make the differential convention
at the point where it is needed.  If $X$ is a real Banach space,
$U\subseteq X$ is open, and $\phi:U\rightarrow\mathbb{R}$, define
\begin{align*}
 D\phi(x)[u]
 =\frac{d}{dt}\bigg|_{t=0}\phi(x+tu),\quad
 D^{2}\phi(x)[u,v]
 =\frac{\partial^{2}}{\partial s\partial t}\bigg|_{s=t=0}
 \phi(x+su+tv).
\end{align*}
For complex noncommutative $L^{q}$-spaces, all derivatives in this
subsection are real Fr\'{e}chet derivatives on the underlying real
Banach space. We need the following semifinite version of the Hessian inequality proved in  \cite{QiuXiong2026}.

\begin{lemma}\cite[Lemma~3.2]{QiuXiong2026}\label{lem:qiu-xiong-hessian-semifinite}
Let $(\mathcal{M},\tau)$ be a semifinite von Neumann algebra. Let $q>2$ and let $p=\frac{q}{q-1}$.  Define $\psi(X)=\|X\|_{L^{q}(\mathcal{M},\tau)}^{p}$ for $X\in L^{q}(\mathcal{M},\tau)$. For $X,Y\in L^{q}(\mathcal{M},\tau)$ with $X\neq0$,
\begin{align*}
 D^{2}\psi(X)[Y,Y]\geq(p-1)D^{2}\psi(X)[iY,iY].
\end{align*}
\end{lemma}

For the passage from the semifinite result to general Haagerup spaces,
we use the following part of the Haagerup--Junge--Xu reduction theorem
\cite[Theorems~2.1 and 3.1]{HaagerupJungeXu2010}.

\begin{lemma}\label{lem:Haagerup-reduction-preliminary}
Let $\mathcal{A}$ be a $\sigma$-finite von Neumann algebra and let
$\varphi$ be a normal faithful state. Set
\begin{align*}
 G_{d}=\bigcup_{n\geq1}2^{-n}\mathbb{Z},\qquad
 \mathcal{R}=\mathcal{A}\rtimes_{\sigma^{\varphi}}G_{d}.
\end{align*}
For $1<q<\infty$, there are finite von Neumann algebras
$\mathcal{R}_{m}\subseteq\mathcal{R}$ with traces $\tau_{m}$,
complex-linear isometries
\begin{align*}
 J_{q}:L^{q}(\mathcal{A})\longrightarrow L^{q}(\mathcal{R}),\qquad
 \iota_{m,q}:L^{q}(\mathcal{R}_{m},\tau_{m})\longrightarrow L^{q}(\mathcal{R}),
\end{align*}
and contractive projections
\begin{align*}
 P_{m}^{(q)}:L^{q}(\mathcal{R})
 \longrightarrow
 \iota_{m,q}\bigl(L^{q}(\mathcal{R}_{m},\tau_{m})\bigr).
\end{align*}
The ranges of $\iota_{m,q}$ are increasing with dense union,
$P_{m}^{(q)}$ fixes every earlier range, and
\begin{align*}
 P_{m}^{(q)}Z\longrightarrow Z
 \qquad\text{in }L^{q}(\mathcal{R})
\end{align*}
for every $Z\in L^{q}(\mathcal{R})$.
\end{lemma}

\begin{lemma}\label{lem:qiu-xiong-hessian}
Let $\mathcal{A}$ be an arbitrary von Neumann algebra. Let $q>2$ and let $p=\frac{q}{q-1}$.  On the Haagerup space $L^{q}(\mathcal{A})$, define $\psi(X)=\|X\|_{L^{q}(\mathcal{A})}^{p}$. For $X,Y\in L^{q}(\mathcal{A})$ with $X\neq0$,
\begin{align*}
 D^{2}\psi(X)[Y,Y]\geq(p-1)D^{2}\psi(X)[iY,iY].
\end{align*}
\end{lemma}
\begin{proof}
Denote $s_{l}(Z)$ and $s_{r}(Z)$ the left and right support projections of $Z\in L^{q}(\mathcal{A})$. Define
\begin{align*}
e=s_{l}(X)\vee s_{r}(X)\vee s_{l}(Y)\vee s_{r}(Y)\quad\mbox{and}\quad \mathcal{A}_{0}=e\mathcal{A}e.
\end{align*}
The left and right support projections of an element of
$L^{q}(\mathcal{A})$ are $\sigma$-finite.  Hence $e$ is
$\sigma$-finite, so is $\mathcal{A}_{0}$, and $X,Y\in
L^{q}(\mathcal{A}_{0})$.  Moreover, the canonical corner inclusion
$I_{e,q}:L^{q}(\mathcal{A}_{0})\rightarrow L^{q}(\mathcal{A})$ is a
complex-linear isometry; see \cite{Terp1982,PisierXu2003}.  Hence
$$\psi_{\mathcal{A}_{0}}=\psi_{\mathcal{A}}\circ I_{e,q}$$
and therefore, by the chain rule,
\begin{align*}
D^{2}\psi_{\mathcal{A}_{0}}(X)[Y,Y]=D^{2}\psi_{\mathcal{A}}(X)[Y,Y],\quad
D^{2}\psi_{\mathcal{A}_{0}}(X)[iY,iY]=D^{2}\psi_{\mathcal{A}}(X)[iY,iY].
\end{align*}
It therefore suffices to prove the assertion for $\mathcal{A}_{0}$. Replacing $\mathcal{A}$ by $\mathcal{A}_{0}$, we may assume from now on that $\mathcal{A}$ is $\sigma$-finite.

Choose a normal faithful state $\varphi$ on $\mathcal{A}$ and let
\begin{align*}
G_{d}=\bigcup_{n\geq 1}2^{-n}\mathbb{Z},\qquad
\mathcal{R}=\mathcal{A}\rtimes_{\sigma^{\varphi}}G_{d}.
\end{align*}
Apply Lemma~\ref{lem:Haagerup-reduction-preliminary} and use the
notation $\mathcal{R}_{m}$, $\tau_{m}$, $J_{q}$,
$\iota_{m,q}$, and $P_{m}^{(q)}$ fixed there.  Its approximation
property gives
\begin{align}\label{eq:xy-convergence}
P_{m}^{(q)}J_{q}(X)\longrightarrow J_{q}(X),\qquad
P_{m}^{(q)}J_{q}(Y)\longrightarrow J_{q}(Y)\quad\mbox{in }L^{q}(\mathcal{R}).
\end{align}
For convenience, denote
\begin{align*}
x=J_{q}(X),\quad y=J_{q}(Y),\quad x_{m}=P_{m}^{(q)}x,\quad y_{m}=P_{m}^{(q)}y.
\end{align*}
Since $X\neq 0$ and $J_{q}$ is isometric, it gives $x\neq 0$. Therefore, $x_{m}\neq 0$ for all sufficiently large $m$. For each such $m$, pick $\widetilde{x}_{m},\widetilde{y}_{m}\in L^{q}(\mathcal{R}_{m},\tau_{m})$ such that
\begin{align*}
x_{m}=\iota_{m,q}(\widetilde{x}_{m}),\qquad y_{m}=\iota_{m,q}(\widetilde{y}_{m}).
\end{align*}

Next, define
\begin{align*}
\psi_{m}(Z)=\lVert Z\rVert_{L^{q}(\mathcal{R}_{m},\tau_{m})}^{p},\quad
\psi_{\mathcal{R}}(Z)=\lVert Z\rVert_{L^{q}(\mathcal{R})}^{p}.
\end{align*}
Observe that $(\mathcal{R}_{m},\tau_{m})$ is a semifinite tracial von Neumann algebra. Applying Lemma \ref{lem:qiu-xiong-hessian-semifinite} to $\widetilde{x}_{m}$ and $\widetilde{y}_{m}$,
\begin{align}
D^{2}\psi_{m}(\widetilde{x}_{m})
[\widetilde{y}_{m},\widetilde{y}_{m}]
\geq
(p-1)
D^{2}\psi_{m}(\widetilde{x}_{m})
[i\widetilde{y}_{m},i\widetilde{y}_{m}].
\label{eq:finite-level-inequality}
\end{align}
Because $\iota_{m,q}$ is a complex-linear isometry, one has $\psi_{m}=\psi_{\mathcal{R}}\circ\iota_{m,q}.$ By chain rule,
\begin{align*}
D^{2}\psi_{m}(\widetilde{x}_{m})[\widetilde{y}_{m},\widetilde{y}_{m}]=D^{2}\psi_{\mathcal{R}}(x_{m})[y_{m},y_{m}],
\quad
D^{2}\psi_{m}(\widetilde{x}_{m})[i\widetilde{y}_{m},i\widetilde{y}_{m}]=D^{2}\psi_{\mathcal{R}}(x_{m})[iy_{m},iy_{m}].
\end{align*}
Consequently, estimate \eqref{eq:finite-level-inequality} becomes
\begin{align*}
D^{2}\psi_{\mathcal{R}}(x_{m})[y_{m},y_{m}]\geq(p-1)D^{2}\psi_{\mathcal{R}}(x_{m})[iy_{m},iy_{m}].
\end{align*}
Since $q>2$ and $p>1$, the function
$Z\mapsto\|Z\|_{q}^{p}$ is twice continuously Fr\'{e}chet differentiable
on $L^{q}(\mathcal{R})\setminus\{0\}$; see
\cite{PotapovSukochevTomskovaZanin2019}.  Consequently,
\begin{align*}
 (Z,W)\longmapsto D^{2}\psi_{\mathcal{R}}(Z)[W,W]
\end{align*}
is continuous on
$(L^{q}(\mathcal{R})\setminus\{0\})\times L^{q}(\mathcal{R})$.
Therefore \eqref{eq:xy-convergence} permits passage to the limit in both
Hessian terms, and we obtain
\begin{align}
D^{2}\psi_{\mathcal{R}}(x)[y,y]\geq(p-1)D^{2}\psi_{\mathcal{R}}(x)[iy,iy].
\label{eq:core-limit-inequality}
\end{align}

Observe that $\psi_{\mathcal{A}}=\psi_{\mathcal{R}}\circ J_{q}$ (see \cite[Theorems 2.1 and 3.1]{HaagerupJungeXu2010}) and $J_{q}(iY)=iJ_{q}(Y)=iy$. By chain rule,
\begin{align*}
D^{2}\psi_{\mathcal{A}}(X)[Y,Y]=D^{2}\psi_{\mathcal{R}}(x)[y,y],\quad
D^{2}\psi_{\mathcal{A}}(X)[iY,iY]=D^{2}\psi_{\mathcal{R}}(x)[iy,iy].
\end{align*}
Combining this with \eqref{eq:core-limit-inequality}, the assertion follows.
\end{proof}

The left Izumi--Hilsum--Terp realization and the Haagerup realization
of $L^{q}(\mathcal{A})$ are canonically complex-linearly isometric.
This isometry preserves multiplication by $i$, the norm, and hence the
first and second real Fr\'{e}chet derivatives of norm functionals; see
\cite{Hilsum1981,Terp1982,Izumi1997,Caspers2013}.  By this canonical
identification,
Lemma~\ref{lem:qiu-xiong-hessian} applies to the left
Izumi--Hilsum--Terp representatives used below.

\subsection{The trace-scaling crossed product on group}
Let $G=N\rtimes_{\alpha}\mathbb{R}$, with the semidirect-product and
Haar-measure conventions fixed in
Section~\ref{subsec:preliminary-solvable-Lie-groups}.
\begin{lemma}\label{lem:group-crossed-product-identification}
Let $N$ be a locally compact group and let $\alpha:\mathbb{R}\longrightarrow\operatorname{Aut}(N)$ be a continuous action. On $G=N\rtimes_{\alpha}\mathbb{R}$, suppose that $d(\alpha_{t}(n))=e^{\omega t}dn.$
Let $\{\beta_{t}\}_{t\in\mathbb{R}}$ be the action on $\operatorname{VN}(N)$ defined by
$$\beta_{t}\bigl(\lambda_{N}(n)\bigr)=\lambda_{N}\bigl(\alpha_{t}(n)\bigr),\quad(n,t)\in N\times\mathbb{R}.$$ 
\begin{enumerate}[\rm(i)]
\item There is a spatial normal $*$-isomorphism $\Theta:\operatorname{VN}(G)\longrightarrow\operatorname{VN}(N)\rtimes_{\beta}\mathbb{R}$ satisfying
$$\Theta\bigl(\lambda_{G}(n,t)\bigr)=\pi_{\beta}\bigl(\lambda_{N}(n)\bigr)\lambda_{t},\qquad (n,t)\in G.$$
In particular,
$$\operatorname{VN}\bigl(N\rtimes_{\alpha}\mathbb{R}\bigr)\simeq\operatorname{VN}(N)\rtimes_{\beta}\mathbb{R}.$$

\item If $f\in C_{c}(G)$ and $f_{t}(n)=f(n,t),$ then
$$\Theta\bigl(\lambda_{G}(f)\bigr)
=\int_{\mathbb{R}}e^{-\omega t}\pi_{\beta}\bigl(\lambda_{N}(f_{t})\bigr)\lambda_{t}\,dt.$$
\end{enumerate}
\end{lemma}
\begin{proof}
For $t\in\mathbb{R}$, define a unitary on $L^{2}(N)$ by
$(U_{t}\xi)(m)=e^{-\frac{\omega t}{2}}\xi(\alpha_{-t}(m)).$
The Jacobian assumption shows that $U_{t}$ is unitary, and a direct
calculation gives
\begin{align*}
 U_{t}\lambda_{N}(n)U_{t}^{*}=\lambda_{N}(\alpha_{t}(n)).
\end{align*}
Thus the formula for $\beta_{t}$ extends uniquely to a spatial normal
$*$-automorphism of $\operatorname{VN}(N)$, and
$t\mapsto\beta_{t}$ is an action.

The group law and Haar measure on $G$ gives the left regular representation on $L^{2}(G,e^{-\omega t}\,dn\,dt)$ by
$$(\lambda_{G}(n_{0},t_{0})\xi)(n,t)=\xi(\alpha_{-t_{0}}(n_{0}^{-1}n),t-t_{0}),\quad(n,t),(n_{0},t_{0})\in N\times\mathbb{R}.$$
Represent the crossed product on $L^{2}(\mathbb{R};L^{2}(N))$ by
$$\bigl(\pi_{\beta}(x)\eta\bigr)(s)=\beta_{-s}(x)\eta(s),\quad(\lambda_{t}\eta)(s)=\eta(s-t),\quad x\in\operatorname{VN}(N),\;\eta\in L^{2}(\mathbb{R};L^{2}(N)).$$
Thus
\begin{align*}
\operatorname{VN}(N)\rtimes_{\beta}\mathbb{R}
=
\left(
\pi_{\beta}(\operatorname{VN}(N))
\cup
\{\lambda_{t}:t\in\mathbb{R}\}
\right)^{\prime\prime}.
\end{align*}
Define $\mathcal{W}: L^{2}(G,e^{-\omega t}\,dn\,dt)\longrightarrow L^{2}(\mathbb{R};L^{2}(N))$ by
$$(\mathcal{W}\xi)(s)(m)=\xi\bigl(\alpha_{s}(m),s\bigr).$$
We first verify that $\mathcal{W}$ is unitary. Indeed,
\begin{align*}
\|\mathcal{W}\xi\|_{L^{2}(\mathbb{R};L^{2}(N))}^{2}
=\int_{\mathbb{R}}\int_{N}\left|\xi\bigl(\alpha_{s}(m),s\bigr)\right|^{2}\,dm\,ds
=\int_{\mathbb{R}}\int_{N}|\xi(n,s)|^{2}e^{-\omega s}\,dn\,ds
=\|\xi\|_{L^{2}(G)}^{2}.
\end{align*}
The inverse of $\mathcal{W}$ is given by
$$(\mathcal{W}^{*}\eta)(n,s)=\eta(s)(\alpha_{-s}(n)).$$
For $n_{0}\in N$, $t_{0}\in\mathbb{R}$, and $\eta\in L^{2}(\mathbb{R};L^{2}(N))$, a direct computation gives that
\begin{align*}
\big(\mathcal{W}\lambda_{G}(n_{0},t_{0})\mathcal{W}^{*}\eta\big)(s)(m)
=\bigl(\lambda_{G}(n_{0},t_{0})\mathcal{W}^{*}\eta\bigr)\bigl(\alpha_{s}(m),s\bigr)
=(\mathcal{W}^{*}\eta)\left(\alpha_{-t_{0}}\bigl(n_{0}^{-1}\alpha_{s}(m)\bigr),s-t_{0}\right).
\end{align*}
Observing that
\begin{align*}
\eta(s-t_{0})\left(\alpha_{-(s-t_{0})}\alpha_{-t_{0}}\bigl(n_{0}^{-1}\alpha_{s}(m)\bigr)\right)
=\eta(s-t_{0})\left(\alpha_{-s}(n_{0})^{-1}m\right),
\end{align*}
we therefore have
\begin{align*}
\big(\mathcal{W}\lambda_{G}(n_{0},t_{0})\mathcal{W}^{*}\eta\big)(s)
=\lambda_{N}\big(\alpha_{-s}(n_{0})\big)\eta(s-t_{0})
=\beta_{-s}\big(\lambda_{N}(n_{0})\big)\eta(s-t_{0})
=\big(\pi_{\beta}(\lambda_{N}(n_{0}))\lambda_{t_{0}}\eta\big)(s).
\end{align*}
Hence
$$\mathcal{W}\lambda_{G}(n_{0},t_{0})\mathcal{W}^{*}=\pi_{\beta}\big(\lambda_{N}(n_{0})\big)\lambda_{t_{0}}.$$
It follows that
\begin{align*}
\mathcal{W}\operatorname{VN}(G)\mathcal{W}^{*}
=\left\{\pi_{\beta}\bigl(\lambda_{N}(n)\bigr)\lambda_{t}:n\in N,\ t\in\mathbb{R}\right\}^{\prime\prime}.
\end{align*}
Taking $t=0$ gives all the generators $\pi_{\beta}\bigl(\lambda_{N}(n)\bigr)$, $n\in N$,
while taking $n=e_{N}$ gives all the generators $\lambda_{t}$. Hence
\begin{align*}
\mathcal{W}\operatorname{VN}(G)\mathcal{W}^{*}
=\left(\pi_{\beta}(\operatorname{VN}(N))\cup\{\lambda_{t}:t\in\mathbb{R}\}\right)^{\prime\prime}
=\operatorname{VN}(N)\rtimes_{\beta}\mathbb{R}.
\end{align*}
Thus $\Theta(x)=\mathcal{W}x\mathcal{W}^{*}$ defines the required spatial normal $*$-isomorphism.

Finally, for $f\in C_{c}(G)$, we have
$\lambda_{G}(f)=\int_{\mathbb{R}}\int_{N}f(n,t)\lambda_{G}(n,t)e^{-\omega t}\,dn\,dt.$
Therefore,
\begin{align*}
\Theta\bigl(\lambda_{G}(f)\bigr)
&=\int_{\mathbb{R}}\int_{N}f(n,t)\pi_{\beta}\bigl(\lambda_{N}(n)\bigr)\lambda_{t}e^{-\omega t}\,dn\,dt\\
&=\int_{\mathbb{R}}e^{-\omega t}\pi_{\beta}\left(\int_{N}f(n,t)\lambda_{N}(n)\,dn\right)\lambda_{t}\,dt\\
&=\int_{\mathbb{R}}e^{-\omega t}\pi_{\beta}\bigl(\lambda_{N}(f_{t})\bigr)\lambda_{t}\,dt.
\end{align*}
The proof is complete.
\end{proof}

The Lemma \ref{lem:group-crossed-product-identification} identifies $\operatorname{VN}\bigl(N\rtimes_{\alpha}\mathbb{R}\bigr)\simeq
\operatorname{VN}(N)\rtimes_{\beta}\mathbb{R}$.
Denote $\Phi$ the dual weight on $\operatorname{VN}(N)\rtimes_{\beta}\mathbb{R}$ associated to trace $\tau_{N}$. For $f\in C_{c}(G)$ and $f_{t}(n)=f(n,t)$, one has
\begin{align}\label{eq:trace-scaling-group-Fourier-element}
 \lambda_{G}(f)=\int_{\mathbb{R}}e^{-\omega t}\pi_{\beta}(\lambda_{N}(f_{t}))\lambda_{t}\,dt.
\end{align}
Indeed, the factor $e^{-\omega t}$ is precisely the Haar density in
\eqref{eq:nonunimodular-haar}.  First let
$f$ belong to the algebraic core
$C_{c}^{\infty}(\mathbb{R})\odot C_{c}^{\infty}(N)$.  Then the
coefficient in \eqref{eq:trace-scaling-group-Fourier-element} belongs
to the coefficient Tomita core, and
Lemma~\ref{lem:trace-scaling-dual-weight} gives
\begin{align*}
 \|\Lambda_{\Phi}(\lambda_{G}(f))\|_{\mathcal{H}_{\Phi}}^{2}
 = \int_{\mathbb{R}} e^{-\omega t}\|\lambda_{N}(f_{t})\|_{L^{2}(\mathcal{M},\tau)}^{2}\,dt
 = \|f\|_{L^{2}(G)}^{2}.
\end{align*}
The algebraic core is dense in $L^{2}(G)$, and its integrated
convolution operators form a left-Hilbert-algebra core for the
Plancherel weight.  Hence the transported dual weight $\Phi$ is the
Plancherel weight of $G$.  The same GNS identity extends by completion
to every $f\in C_{c}(G)$.

The standard Fourier--Lebesgue and left Izumi--Hilsum--Terp
identifications in Sections~\ref{subsec:preliminary-group-Fourier-transform}
and~\ref{subsec:preliminary-dual-weights} now give, on the
Fourier--Tomita cores,
\begin{align*}
 \mathcal{F}_{G,p}f\longleftrightarrow j_{q}^{\Phi}(\lambda_{G}(f)),\qquad
 \mathcal{F}_{N,p}h \longleftrightarrow j_{q}^{\tau_{N}}(\lambda_{N}(h)).
\end{align*}
For general $L^{p}$ inputs, these expressions denote the continuous
interpolation extensions.  In particular, the source core element is
$j_{q}^{\tau_{N}}(\lambda_{N}(h))$, with Fourier-side index $q$.

\begin{corollary}\label{cor:coarse-trace-scaling-extension}
Under assumptions in Lemma \ref{lem:group-crossed-product-identification}, let $1\leq p\leq2$ and $q=\frac{p}{p-1}$, with $q=\infty$ when
$p=1$.
For every $f\in L^{p}(G)$,
\begin{align*}
\|\mathcal{F}_{G,p}f\|_{\mathcal{F}L^{q}(G)}\leq H_{p}(N)\|f\|_{L^{p}(G)}.
\end{align*}
In particular, the same estimate holds with $H_{p}(N)$ replaced by
$1$.
\end{corollary}

\begin{proof}
Let $f$ first belong to the common Fourier--Tomita core.  At the
$p=1$ endpoint, \eqref{eq:trace-scaling-group-Fourier-element} is the
identity
\begin{align*}
 \mathcal{A}_{\omega}\bigl(t\longmapsto
 j_{\infty}^{\tau_{N}}(\lambda_{N}(f_{t}))\bigr)
 =\mathcal{A}_{\omega}\bigl(t\longmapsto\lambda_{N}(f_{t})\bigr)
 =\lambda_{G}(f)
 =j_{\infty}^{\Phi}(\lambda_{G}(f)).
\end{align*}
At the $p=2$ endpoint, Lemma
\ref{lem:trace-scaling-dual-weight} gives, under the GNS
identifications,
\begin{align*}
 \mathcal{A}_{\omega}\bigl(t\longmapsto
 j_{2}^{\tau_{N}}(\lambda_{N}(f_{t}))\bigr)
 =j_{2}^{\Phi}(\lambda_{G}(f)).
\end{align*}
Compatibility and uniqueness of complex interpolation therefore give,
for $1<p<2$,
\begin{align*}
 \mathcal{A}_{\omega}\bigl(t\longmapsto
 j_{q}^{\tau_{N}}(\lambda_{N}(f_{t}))\bigr)
 =j_{q}^{\Phi}\left(
 \int_{\mathbb{R}}e^{-\omega t}
 \pi_{\beta}(\lambda_{N}(f_{t}))\lambda_{t}\,dt
 \right)
 =j_{q}^{\Phi}(\lambda_{G}(f)).
\end{align*}
Taking norms in the concrete realizations gives
\begin{align*}
 \|\mathcal{F}_{G,p}f\|_{\mathcal{F}L^{q}(G)}
=\|j_{q}^{\Phi}(\lambda_{G}(f))\|_{L^{q}(\operatorname{VN}(N)\rtimes_{\beta}\mathbb{R})_{\mathrm{left}}},
\quad \|\mathcal{F}_{N,p}f_{t}\|_{\mathcal{F}L^{q}(N)}
 =\|j_{q}^{\tau_{N}}(\lambda_{N}(f_{t}))\|_{L^{q}(\mathcal{M},\tau_{N})}.
\end{align*}
Equivalently, after transporting the assembly map through the canonical
Fourier--Lebesgue identifications,
\begin{align*}
 \mathcal{A}_{\omega}\bigl(t\longmapsto\mathcal{F}_{N,p}f_{t}\bigr)
 =\mathcal{F}_{G,p}f
 \quad\text{in}\quad\mathcal{F}L^{q}(G).
\end{align*}
This is an identity for the interpolated assembly map; it does not insert
the generally unbounded $L^{q}$-element $\mathcal{F}_{N,p}f_{t}$ into
$\pi_{\beta}$ pointwise.  Applying Lemma
\ref{lem:weighted-crossed-product-assembly},
\begin{align*}
\|\mathcal{F}_{G,p}f\|_{\mathcal{F}L^{q}(G)}
\leq\left(\int_{\mathbb{R}}e^{-\omega t}
\|\mathcal{F}_{N,p}f_{t}\|_{\mathcal{F}L^{q}(N)}^{p}\,dt\right)^{\frac{1}{p}}
\leq H_{p}(N)\|f\|_{L^{p}(G)}.
\end{align*}  
Finally, density extends the assertion to every $f\in L^{p}(G)$. The proof is complete.
\end{proof}

\subsection{Uniform window bounds and density}
Denote $\Phi$ the dual weight on $\operatorname{VN}(N)\rtimes_{\beta}\mathbb{R}$ associated to trace $\tau_{N}$.
For $h\in C_{c}^{\infty}(N)$ and $\varepsilon>0$, define
\begin{align}
 h_{\varepsilon}
 =\left(\frac{\varepsilon}{\pi}\right)^{\frac{1}{2}}\int_{\mathbb{R}}e^{-\varepsilon r^{2}}V_{r}h\,dr,
 \label{eq:flow-analytic-regularization}
\end{align}
where $V_{r}$ is given by setting
\begin{align*}
 (V_{t}h)(n)=e^{-\omega t}h(\alpha_{-t}(n)).
\end{align*}
A change of variables gives that
\begin{align*}
 \|V_{t}h\|_{L^{v}(N)}
 =e^{-\omega t\left(1-\frac{1}{v}\right)}\|h\|_{L^{v}(N)}.
\end{align*}
The entire extension $V_{z}h_{\varepsilon}
 =\left(\frac{\varepsilon}{\pi}\right)^{\frac{1}{2}}\int_{\mathbb{R}}e^{-\varepsilon(r-z)^{2}}V_{r}h\,dr$, $z\in\mathbb{C}$, satisfies
\begin{align*}
 \|V_{z}h_{\varepsilon}\|_{L^{v}(N)}
 \leq\exp\left(\varepsilon(\operatorname{Im}z)^{2}-c_{v}\operatorname{Re}z+\frac{c_{v}^{2}}{4\varepsilon}\right)
 \|h\|_{L^{v}(N)},\quad 1\leq v<\infty,
\end{align*}
where $c_{v}=\omega\left(1-\frac{1}{v}\right)$. 

Let $\mathfrak{A}_{N}$ be the algebra generated under convolution and involution by the functions in
\eqref{eq:flow-analytic-regularization}.  Since
\begin{align*}
 \lambda_{N}(V_{r}h)
 =
 \beta_{r}(\lambda_{N}(h)),
\end{align*}
$\lambda_{N}(\mathfrak{A}_{N})$ is a $\beta$-entire Tomita algebra. Moreover, the normalized Gaussians form an approximate identity for the strongly continuous groups $\{V_{r}\}_{r\in\mathbb{R}}$ on the relevant $L^{v}$-spaces.
Hence $\mathfrak{A}_{N}$ is dense in $L^{v}(N)$, $1\leq v<\infty$. This is the usual Gaussian construction of analytic elements in a left Hilbert algebra; see, for example, \cite[Chapter~VIII]{Takesaki2003} and \cite{Terp1982}.

\begin{lemma}\label{lem:analytic-flow-core}
Under the assumptions of Theorem \ref{thm:one-dimensional-extension},
let $G=N\rtimes_{\alpha}\mathbb{R}$, let $1\leq p\leq2$, and let
$q=\frac{p}{p-1}$, with $q=\infty$ when $p=1$. We have
\begin{enumerate}[\rm(i)]
\item\label{lem:analytic-flow-core-density} The algebraic tensor product
$\mathfrak{A}_{N}\odot C_{c}^{\infty}(\mathbb{R})$ is dense in
$L^{v}(G)$ for $1\leq v<\infty$.
\item\label{lem:analytic-flow-core-function} If $f\in \mathfrak{A}_{N}\odot C_{c}^{\infty}(\mathbb{R})$ and
$f_{a,b}(n,t):=f(n,t)e^{-at^{2}+bt}$ for $a\geq0,$ $b\in\mathbb{R},$
then $$(a,b)\longmapsto \mathcal{F}_{G,p}(f_{a,b})$$ is twice continuously differentiable as a map into
$\mathcal{F}L^{q}(G)$ on $(0,\infty)\times\mathbb{R}$, with continuous one-sided derivatives at $a=0$, and
\begin{align*}
  \partial_{b}\mathcal{F}_{G,p}(f_{a,b})=\mathcal{F}_{G,p}(t f_{a,b}),\quad
  \partial_{b}^{2}\mathcal{F}_{G,p}(f_{a,b})=\mathcal{F}_{G,p}(t^{2} f_{a,b}),\quad
  \partial_{a}\mathcal{F}_{G,p}(f_{a,b}) =-\mathcal{F}_{G,p}(t^{2} f_{a,b}).
\end{align*}
\item\label{replace to other} All the corresponding identities remain valid after applying a
quotient modulation
$f(n,t)\longmapsto e^{i\rho t}f(n,t)$ for any
$\rho\in\mathbb{R}$.  Moreover, for fixed $(a,b)$, the map
\begin{align*}
 r\longmapsto\mathcal{F}_{G,p}(e^{irt}f_{a,b})
\end{align*}
is twice continuously differentiable, with first and second
derivatives at $r=0$ equal to
$i\mathcal{F}_{G,p}(tf_{a,b})$ and
$-\mathcal{F}_{G,p}(t^{2}f_{a,b})$, respectively.
\end{enumerate}
\end{lemma}
\begin{proof}
The space $\mathfrak{A}_{N}$ is dense in $L^{v}(N)$ for
$1\leq v<\infty$, and
$C_{c}^{\infty}(\mathbb{R})$ is dense in
$L^{v}(\mathbb{R},e^{-\omega t}\,dt)$.  The Bochner-space
tensor-product density theorem
\cite{BerghLofstrom1976,PisierXu2003} therefore gives assertion
\eqref{lem:analytic-flow-core-density}.

For $f$ in the displayed algebraic tensor product, its support in the
$t$-variable is compact.  Differentiation of the scalar multiplier
$e^{-at^{2}+bt}$, followed by dominated convergence on this compact
support, shows that $(a,b)\mapsto f_{a,b}$ is twice continuously
differentiable in $L^{p}(G)$, with continuous one-sided derivatives at
$a=0$, and
\begin{align*}
 \partial_{b}f_{a,b}=tf_{a,b},\qquad
 \partial_{b}^{2}f_{a,b}=t^{2}f_{a,b},\qquad
 \partial_{a}f_{a,b}=-t^{2}f_{a,b}.
\end{align*}
Indeed, $\lambda_{N}(\mathfrak{A}_{N})$ is a $\beta$-entire Tomita
algebra and the coefficient functions have compact smooth support in the
$t$-variable.  Hence they lie in the common Fourier--Tomita core, and
the bounded linear map $\mathcal{F}_{G,p}$ carries these identities into
$\mathcal{F}L^{q}(G)$.  Under the canonical concrete realization, these
are precisely the corresponding derivative identities for the
$j_{q}^{\Phi}$-representatives in
$L^{q}(\operatorname{VN}(N)\rtimes_{\beta}\mathbb{R})_{\mathrm{left}}$.  This
proves \eqref{lem:analytic-flow-core-function}.  Multiplication by
$e^{i\rho t}$ preserves the compact $t$-support and commutes with all
three differentiations.  The same compact-support dominated-convergence
argument in the modulation parameter gives the two displayed
$r$-derivatives, proving \eqref{replace to other}.
\end{proof}

The next quotient modulation invariance is the bridge between the imaginary part of the Hessian and differentiation in the quotient variable.

\begin{lemma}\label{lem:quotient-modulation}
Under the assumptions in Theorem \ref{thm:one-dimensional-extension},
let $G=N\rtimes_{\alpha}\mathbb{R}$, let $1\leq p\leq2$, and let
$q=\frac{p}{p-1}$, with $q=\infty$ when $p=1$. For
$\rho\in\mathbb{R}$, define
\begin{align*}
  (M_{\rho} f)(n,t)=e^{i\rho t}f(n,t),\quad n\in N,\;t\in\mathbb{R}.
\end{align*}
Then
$$\|\mathcal{F}_{G,p}(M_{\rho} f)\|_{\mathcal{F}L^{q}(G)}
  =\|\mathcal{F}_{G,p}f\|_{\mathcal{F}L^{q}(G)},\quad f\in L^{p}(G).$$
\end{lemma}
\begin{proof}
By Lemma \ref{lem:group-crossed-product-identification}, identify $\operatorname{VN}(G)\simeq\operatorname{VN}(N)\rtimes_{\beta}\mathbb{R}$. The
dual action is given on generators by
\begin{align*}
  \widehat{\beta}_{\rho}(\pi_{\beta}(x))=\pi_{\beta}(x),
  \qquad
  \widehat{\beta}_{\rho}(\lambda_{t})=e^{i\rho t}\lambda_{t}.
\end{align*}
The dual weight is invariant under this action:
\begin{align*}
 \Phi\circ\widehat{\beta}_{\rho}=\Phi.
\end{align*}
It is therefore isometric both on $\operatorname{VN}(N)\rtimes_{\beta}\mathbb{R}$ and on
$\mathcal{H}_{\Phi}$, because for $Y\in\mathfrak{n}_{\Phi}$,
\begin{align*}
 \|\Lambda_{\Phi}(\widehat{\beta}_{\rho}(Y))\|_{\mathcal{H}_{\Phi}}^{2}
 =\Phi(\widehat{\beta}_{\rho}(Y^{*}Y))
 =\Phi(Y^{*}Y).
\end{align*}
Interpolation of the compatible endpoint actions consequently gives an
isometry on $L^{q}(\operatorname{VN}(N)\rtimes_{\beta}\mathbb{R})_{\mathrm{left}}$, satisfying on the common
core
\begin{align*}
 \widehat{\beta}_{\rho}\bigl(j_{q}^{\Phi}(Y)\bigr)
 =j_{q}^{\Phi}(\widehat{\beta}_{\rho}(Y)).
\end{align*}
The coefficient formula \eqref{eq:trace-scaling-group-Fourier-element}
now gives
\begin{align*}
 j_{q}^{\Phi}(\lambda_{G}(M_{\rho}f))
 =\widehat{\beta}_{\rho}\bigl(j_{q}^{\Phi}(\lambda_{G}(f))\bigr)
\end{align*}
on the common core.  Through the canonical isometry with
$\mathcal{F}L^{q}(G)$, density gives the asserted norm equality for
every $f\in L^{p}(G)$.
\end{proof}

For the remainder of the flow and freezing argument, assume
$1<p<2$ and $q=\frac{p}{p-1}$.  Let
$c_{\rho}=\left(\frac{p\rho}{\pi}\right)^{\frac{1}{2p}}$ be chosen so that
$ c_{\rho}^{p}\int_{\mathbb{R}}e^{-p\rho(t-s)^{2}}\,ds=1$ for $\rho>0$. The following result provides core monotonicity and endpoint estimate for flow.
\begin{lemma}\label{lem:flow-core-monotonicity}
Under assumptions in Theorem \ref{thm:one-dimensional-extension}, let $G=N\rtimes_{\alpha}\mathbb{R}$. For $\rho>0$ and $f\in\mathfrak{A}_{N}\odot C_{c}^{\infty}(\mathbb{R})$, define
$$\mathcal{L}_{\rho}(f)
  =c_{\rho}^{p}\int_{\mathbb{R}} \|\mathcal{F}_{G,p}(W_{\rho,s}f)\|_{\mathcal{F}L^{q}(G)}^{p}\,ds,$$
where $$c_{\rho}=\left(\frac{p\rho}{\pi}\right)^{\frac{1}{2p}},\quad W_{\rho,s}f(n,t)=e^{-\rho(t-s)^{2}}f(n,t),\quad s\in\mathbb{R}.$$
We have
\begin{enumerate}[\rm(i)]
    \item\label{monotonicity property} monotonicity property: for $0<\rho_{1}<\rho_{2}$,
    $$\mathcal{L}_{\rho_{1}}(f)\leq\mathcal{L}_{\rho_{2}}(f);$$
    \item\label{limitation property} limitation property: 
    $$\lim_{\rho\downarrow0}\mathcal{L}_{\rho}(f) =\|\mathcal{F}_{G,p}f\|_{\mathcal{F}L^{q}(G)}^{p}.$$
\end{enumerate}
\end{lemma}

\begin{proof}
The assertion is immediate for $f=0$, so assume $f\neq0$.  Let $f_{a,b}$ be given as in Lemma \ref{lem:analytic-flow-core} \eqref{lem:analytic-flow-core-function}. Denote $Y(a,b)=j_{q}^{\Phi}(\lambda_{G}(f_{a,b}))$ and
\begin{align*}
V(a,b)=\|Y(a,b)\|_{L^{q}(\operatorname{VN}(N)\rtimes_{\beta}\mathbb{R})_{\mathrm{left}}}^{p}
  =\|\mathcal{F}_{G,p}(f_{a,b})\|_{\mathcal{F}L^{q}(G)}^{p}.
\end{align*}
The element $Y(a,b)$ is nonzero.  Indeed, all the functions involved
belong to the Fourier--Tomita core, and
\begin{align*}
 Y(a,b)=0
 \Longrightarrow \lambda_{G}(f_{a,b})=0
 \Longrightarrow f_{a,b}=0
 \Longrightarrow f=0,
\end{align*}
where the first implication uses the injectivity of $j_{q}^{\Phi}$ on the
common core and the second follows from the injectivity of the Plancherel
GNS map.  This contradicts the choice of $f$.

For $\psi(T)=\|T\|_{L^{q}(\operatorname{VN}(N)\rtimes_{\beta}\mathbb{R})_{\mathrm{left}}}^{p}$, the
chain rule and Lemma \ref{lem:analytic-flow-core}
\eqref{lem:analytic-flow-core-function} give
\begin{align}  \label{eq:flow-second-b-chain-rule}
  \partial_{b}^{2}V(a,b)
  =D^{2}\psi(Y(a,b))[X(a,b),X(a,b)]+D\psi(Y(a,b))[Z(a,b)]
\end{align}
where
\begin{align*}
  X(a,b)=j_{q}^{\Phi}(\lambda_{G}(t f_{a,b}))\quad\mbox{and}\quad
  Z(a,b)=j_{q}^{\Phi}(\lambda_{G}(t^{2} f_{a,b})).
\end{align*}
On the other hand, Lemma~\ref{lem:quotient-modulation} says that
\begin{align*}
  r\longmapsto
  \psi\bigl(j_{q}^{\Phi}(\lambda_{G}(e^{irt}f_{a,b}))\bigr)
  =\psi\bigl(j_{q}^{\Phi}(\lambda_{G}(f_{a,b}))\bigr).
\end{align*}
Twice differentiating at $r=0$ gives
\begin{align} \label{eq:flow-modulation-second-derivative}
  D^{2}\psi(Y(a,b))[iX(a,b),iX(a,b)] =D\psi(Y(a,b))[Z(a,b)].
\end{align}
Combining \eqref{eq:flow-second-b-chain-rule}, \eqref{eq:flow-modulation-second-derivative} and Lemma \ref{lem:qiu-xiong-hessian}
\begin{align*}
  \partial_{b}^{2}V(a,b)&\geq(p-1)D^{2}\psi(Y(a,b))[iX(a,b),iX(a,b)]+D\psi(Y(a,b))[Z(a,b)]\\
  &=pD\psi(Y(a,b))[Z(a,b)].
\end{align*}
Since $\partial_{a}Y(a,b)=-Z(a,b)$, we have 
\begin{align}\label{formula:pb2+pa V}
  \bigl(\partial_{b}^{2}+p\partial_{a}\bigr)V(a,b)\geq0.
\end{align}

If the $t$-support of $f$ is contained in $[-t_{0},t_{0}]$, the
contractivity of $\mathcal{F}_{G,p}$ and the derivative formulas above
give, locally uniformly for $a\geq0$,
\begin{align*}
 |V(a,b)|+|\partial_{a}V(a,b)|+|\partial_{b}^{2}V(a,b)|
 \leq C e^{C|b|}
\end{align*}
for a constant $C$ depending on the chosen compact $a$-interval and on
$f$.  Thus all heat-kernel integrals below are finite, and dominated
convergence justifies differentiation under those integrals.

Let $H_{r}=e^{r\partial_{b}^{2}}$ be the heat semigroup on $\mathbb{R}$. Fix $0\leq a<\rho$, $b\in\mathbb{R}$, and $\varepsilon>0$. For $0\leq r\leq \rho-a$, one obtains
\begin{align*}
 \frac{d}{dr}H_{\frac{r+\varepsilon}{p}}\Big(V(a+r,\cdot)\Big)(b)
 =\frac{1}{p} H_{\frac{r+\varepsilon}{p}}\Big((\partial_{b}^{2}+p\partial_{a})V(a+r,\cdot)\Big)(b)\geq0
\end{align*}
where we use \eqref{formula:pb2+pa V} for the inequality. Hence, $r\mapsto H_{\frac{r+\varepsilon}{p}}\Big(V(a+r,\cdot)\Big)(b)$ is non-decreasing and we therefore have
\begin{align*}
 H_{\frac{\varepsilon}{p}}\Big(V(a,\cdot)\Big)(b)\leq H_{\frac{\rho-a+\varepsilon}{p}}\Big(V(\rho,\cdot)\Big)(b).
\end{align*}
Letting $\varepsilon\downarrow0$ gives
\begin{align}\label{eq:flow-heat-comparison}
  V(a,b)\leq H_{\frac{\rho-a}{p}}\Big(V(\rho,\cdot)\Big)(b).
\end{align}
Observe that
\begin{align*}
  \mathcal{L}_{\rho}(f)=H_{\frac{\rho}{p}}\Big(V(\rho,\cdot)\Big)(0).
\end{align*}
Combining this with \eqref{eq:flow-heat-comparison},
\begin{align*}
  \mathcal{L}_{\rho_{1}}(f)
  \leq H_{\frac{\rho_{1}}{p}}H_{\frac{\rho_{2}-\rho_{1}}{p}}\Big(V(\rho_{2},\cdot)\Big)(0)
  =\mathcal{L}_{\rho_{2}}(f),
\end{align*}
which gives the monotonicity property \eqref{monotonicity property}.

Finally, if $0<\rho\leq \rho_{0}$ and the $t$-support of $f$ is contained in $[-t_{0},t_{0}]$, then the Hausdorff--Young inequality in Corollary \ref{cor:coarse-trace-scaling-extension} and definition of $f_{a,b}$ give that
\begin{align*}
  e^{-t^{2}}V\big(\rho,2\sqrt{\frac{\rho}{p}}t\big)\leq e^{-t^{2}}\|f_{\rho,2\sqrt{\frac{\rho}{p}}t}\|_{L^{p}(G)}^{p}
  \leq e^{-t^{2}+2t_{0}\sqrt{p\rho_{0}}|t|}\|f\|_{L^{p}(G)}^{p}.
\end{align*}
Thus, the dominated convergence theorem implies that
\begin{align*}
\lim_{\rho\downarrow0}\mathcal{L}_{\rho}(f)
  =\frac{1}{\sqrt{\pi}}\int_{\mathbb{R}}e^{-t^{2}}\lim_{\rho\downarrow0} V\big(\rho,2\sqrt{\frac{\rho}{p}}t\big)\,dt
  =V(0,0)=\|\mathcal{F}_{G,p}f\|_{\mathcal{F}L^{q}(G)}^{p}.
\end{align*}
This proves limitation property \eqref{limitation property} and the proof is complete.
\end{proof}

The normalization for $W_{\rho,s}$ makes the family of window operators uniformly contractive.  This is what permits the flow obtained on the density set to extend without loss of constant.

\begin{lemma}\label{lem:flow-window-density}
Under assumptions in Theorem \ref{thm:one-dimensional-extension}, let $G=N\rtimes_{\alpha}\mathbb{R}$ and let $\mathcal{L}_{\rho}(f)$ be defined as in Lemma \ref{lem:flow-core-monotonicity}.  For $f\in L^{p}(G)$, the function $\rho\mapsto\mathcal{L}_{\rho}(f)$ is nondecreasing and
$$\lim_{\rho\downarrow0}\mathcal{L}_{\rho}(f)=\|\mathcal{F}_{G,p}f\|_{\mathcal{F}L^{q}(G)}^{p}.$$
\end{lemma}
\begin{proof}
For every $\rho>0$, let $c_{\rho}$ be given as in Lemma \ref{lem:flow-core-monotonicity} and define  the map
\begin{align*}
  T_{\rho}f(s)=c_{\rho}\mathcal{F}_{G,p}(W_{\rho,s}f),\quad f\in\mathfrak{A}_{N}\odot C_{c}^{\infty}(\mathbb{R}).
\end{align*}
The map $s\mapsto W_{\rho,s}f$ is strongly measurable in $L^{p}(G)$; so is $T_{\rho}f$.  By the Hausdorff--Young inequality in Corollary \ref{cor:coarse-trace-scaling-extension} and the definition of $c_{\rho}$,
\begin{align*}
  \|T_{\rho}f\|_{L^{p}\bigl(\mathbb{R};\mathcal{F}L^{q}(G)\bigr)}^{p}
  \leq c_{\rho}^{p}\int_{\mathbb{R}}\int_{\mathbb{R}}\int_{N} e^{-\omega t}e^{-p\rho(t-s)^{2}}|f(n,t)|^{p}\,dn\,dt\,ds
  =\|f\|_{L^{p}(G)}^{p}.
\end{align*}
Thus, $T_{\rho}$ has a contractive extension from $L^{p}(G)$ to
$L^{p}\bigl(\mathbb{R};\mathcal{F}L^{q}(G)\bigr)$. Moreover, for
$f,h\in L^{p}(G)$,
\begin{align*}
  \Big|\mathcal{L}_{\rho}(f)^{\frac{1}{p}}-\mathcal{L}_{\rho}(h)^{\frac{1}{p}}\Big|
  =\Big|\|T_{\rho}f\|_{L^{p}\bigl(\mathbb{R};\mathcal{F}L^{q}(G)\bigr)}
  -\|T_{\rho}h\|_{L^{p}\bigl(\mathbb{R};\mathcal{F}L^{q}(G)\bigr)}\Big|
  \leq\|f-h\|_{L^{p}(G)}
\end{align*}
uniformly in $\rho>0$.  Choose
$f_{n}\in\mathfrak{A}_{N}\odot C_{c}^{\infty}(\mathbb{R})$ with
$f_{n}\to f$ in $L^{p}(G)$.  Then
\begin{align*}
 \left|\mathcal{L}_{\rho}(f)^{\frac{1}{p}}
 -\|\mathcal{F}_{G,p}f\|_{\mathcal{F}L^{q}(G)}\right|
 &\leq \|f-f_{n}\|_{p}
 +\left|\mathcal{L}_{\rho}(f_{n})^{\frac{1}{p}}
 -\|\mathcal{F}_{G,p}f_{n}\|_{\mathcal{F}L^{q}(G)}\right|\\
 &\qquad+\|\mathcal{F}_{G,p}(f_{n}-f)\|_{\mathcal{F}L^{q}(G)}.
\end{align*}
The first and third terms are bounded by $\|f-f_{n}\|_{p}$,
uniformly in $\rho$.  Lemma~\ref{lem:flow-core-monotonicity}, followed
by $n\to\infty$, proves the small-window limit.  The same uniform
estimate transfers monotonicity from the dense core to all of
$L^{p}(G)$.
\end{proof}

\subsection{Weighted freezing and the large-window endpoint}

 For $f\in L^{p}(N\rtimes_{\alpha}\mathbb{R}),$
choose a strongly
measurable representative $t\mapsto f_{t}$, where $f_{t}(n)=f(n,t)$. By Lemma \ref{lem:group-crossed-product-identification}, identify $L^{p}(N\rtimes_{\alpha}\mathbb{R})
 =L^{p}(\mathbb{R},e^{-\omega t}\,dt;L^{p}(N))$.
For $r\in\mathbb{R}$, define
\begin{align*}
  D_{f}(r)=\int_{\mathbb{R}}e^{-\omega t}
  \|f_{t}-f_{t-r}\|_{L^{p}(N)}^{p}\,dt.
\end{align*}
The translation group is strongly continuous on this weighted Bochner
space, and a change of variables gives
\begin{align*}
  \lim_{r\to0}D_{f}(r)=0,\qquad  D_{f}(r)\leq2^{p-1}(1+e^{-\omega r})\|f\|_{L^{p}(G)}^{p}.
\end{align*}
Splitting the integral into a neighborhood of zero and its complement,
and using the Gaussian to dominate the possible exponential growth,
therefore gives
\begin{align}\label{formula:integral convergence}
\rho^{\frac{1}{2}}\int_{\mathbb{R}}e^{-p\rho r^{2}}D_{f}(r)\,dr\longrightarrow0,\quad\mbox{as}\quad\rho\to\infty.
\end{align}

\begin{lemma}\label{lem:trace-preserving-freezing}
Under assumptions in Theorem \ref{thm:one-dimensional-extension}, let $G=N\rtimes_{\alpha}\mathbb{R}$ and let $\mathcal{L}_{\rho}(f)$ be defined as in Lemma \ref{lem:flow-core-monotonicity}. For $f\in L^{p}(G)$, we have
$$\limsup_{\rho\to\infty}\mathcal{L}_{\rho}(f)^{\frac{1}{p}}\leq B_{p} H_{p}(N)\|f\|_{L^{p}(G)}.$$
\end{lemma}
\begin{proof}
Denote $g_{\rho,s}(t)=e^{-\rho(t-s)^{2}}$. For $s\in\mathbb{R}$, define
\begin{align*}
  R_{\rho}(s)=c_{\rho}\|\mathcal{F}_{G,p}(W_{\rho,s}f)\|_{\mathcal{F}L^{q}(G)},\qquad
  S_{\rho}(s)=c_{\rho}\|\mathcal{F}_{G,p}(f_{s}\otimes g_{\rho,s})\|_{\mathcal{F}L^{q}(G)},
\end{align*}
where $f_{t}(n)=f(n,t)$, $n\in N$, $t\in\mathbb{R}$.

The map $s\mapsto W_{\rho,s}f$ is continuous in $L^{p}(G)$, so
$R_{\rho}$ is measurable.  The map
$s\mapsto\mathcal{F}_{N,p}f_{s}$ is strongly measurable in
$\mathcal{F}L^{q}(N)$; equivalently, its concrete representative is
strongly measurable in $L^{q}(\operatorname{VN}(N),\tau_{N})$.  Moreover,
$s\longmapsto
 \widehat{e^{-\frac{\omega(\cdot)}{p}}g_{\rho,s}(\cdot)}$
is continuous in $L^{q}(\mathbb{R})$, where $\widehat{\cdot}$ denotes the scalar Fourier transform defined in \eqref{eq:scalar-Fourier-transform}.

For a Fourier--Tomita core element $h$ on $N$, the identity
$\frac{1}{p}+\frac{1}{q}=1$ and
\eqref{eq:trace-scaling-group-Fourier-element} give
\begin{align*}
 \lambda_{G}(h\otimes g_{\rho,s})
 &=\int_{\mathbb{R}}e^{-\omega t}g_{\rho,s}(t)
 \pi_{\beta}(\lambda_{N}(h))\lambda_{t}\,dt\\
 &=\int_{\mathbb{R}}e^{-\frac{\omega t}{q}}
 \Big(e^{-\frac{\omega t}{p}}g_{\rho,s}(t)\Big)
 \pi_{\beta}(\lambda_{N}(h))\lambda_{t}\,dt.
\end{align*}
The scalar function $e^{-\frac{\omega(\cdot)}{p}}g_{\rho,s}$ need not have
compact support.  Choose $b_{k}\in C_{c}^{\infty}(\mathbb{R})$ such that
\begin{align*}
 e^{-\frac{\omega(\cdot)}{q}}b_{k}
 \longrightarrow e^{-\omega(\cdot)}g_{\rho,s}
 \quad\text{in }L^{1}(\mathbb{R})\quad\mbox{and}\quad
 \widehat{b_{k}}\longrightarrow
 \widehat{e^{-\frac{\omega(\cdot)}{p}}g_{\rho,s}}
 \quad\text{in }L^{q}(\mathbb{R}).
\end{align*}
Such a sequence is obtained by smooth cutoffs of this Gaussian.  The
first, weighted $L^{1}$ convergence identifies the operator-valued
coefficient integral,
while the $L^{q}$ convergence and continuity in Theorem
\ref{thm:weighted-crossed-product-kss} identify its interpolated image.
Under the core identifications, Theorem
\ref{thm:weighted-crossed-product-kss} therefore sends
\begin{align*}
 \left(
 j_{q}^{\tau_{N}}(\lambda_{N}(h)),
 \widehat{e^{-\frac{\omega(\cdot)}{p}}g_{\rho,s}(\cdot)}
 \right)
\end{align*}
to the concrete representative
\begin{align*}
 j_{q}^{\Phi}(\lambda_{G}(h\otimes g_{\rho,s}))
 \in L^{q}(\operatorname{VN}(N)\rtimes_{\beta}\mathbb{R})_{\mathrm{left}}.
\end{align*}
The two Fourier--Lebesgue norms are therefore represented by
$$\|\mathcal{F}_{G,p}(h\otimes g_{\rho,s})\|_{\mathcal{F}L^{q}(G)}
 =\|j_{q}^{\Phi}(\lambda_{G}(h\otimes g_{\rho,s}))\|_{L^{q}(\operatorname{VN}(N)\rtimes_{\beta}\mathbb{R})_{\mathrm{left}}}$$
 and
$$ \|\mathcal{F}_{N,p}h\|_{\mathcal{F}L^{q}(N)}
 =\|j_{q}^{\tau_{N}}(\lambda_{N}(h))\|_{L^{q}(\operatorname{VN}(N),\tau_{N})}.$$
The core is dense in
$L^{p}(N)$, and both the Kato--Seiler--Simon bilinear map and the map
$h\mapsto\mathcal{F}_{G,p}(h\otimes g_{\rho,s})$ are continuous.
Consequently the identity extends to $h=f_{s}$ for almost every $s$;
no expression $\lambda_{N}(f_{s})$ is used for a general
$L^{p}$-slice.  The strong measurability of the first input, the
continuity of the second input, and continuity of the Kato--Seiler--Simon
bilinear map now also prove that $S_{\rho}$ is measurable.  Hence
\begin{align*}
\|\mathcal{F}_{G,p}(f_{s}\otimes g_{\rho,s})\|_{\mathcal{F}L^{q}(G)}
&\leq\|\mathcal{F}_{N,p}f_{s}\|_{\mathcal{F}L^{q}(N)}\|\widehat{e^{-\frac{\omega(\cdot)}{p}}g_{\rho,s}(\cdot)}\|_{L^{q}(\mathbb{R})}\\
&\leq H_{p}(N)
\|f_{s}\|_{L^{p}(N)}\|\widehat{e^{-\frac{\omega(\cdot)}{p}}g_{\rho,s}(\cdot)}\|_{L^{q}(\mathbb{R})}.
\end{align*}
Since $c_{\rho}\left\|\widehat{e^{-\rho(\cdot)^{2}}}\right\|_{L^{q}(\mathbb{R})} =p^{\frac{1}{2p}}q^{-\frac{1}{2q}}=B_{p}$,
we have
\begin{align*}
 S_{\rho}(s)\leq \exp\left(\frac{\omega^{2}}{4p^{2}\rho}\right)e^{-\frac{\omega s}{p}}B_{p}H_{p}(N)\|f_{s}\|_{L^{p}(N)}.
\end{align*}
Moreover,
\begin{align}\label{eq:unified-frozen-field-bound}
 \|S_{\rho}\|_{L^{p}(\mathbb{R})}\leq \exp\left(\frac{\omega^{2}}{4p^{2}\rho}\right) B_{p}H_{p}(N)\|f\|_{L^{p}(G)}.
\end{align}

Finally, by the Hausdorff--Young inequality in Corollary \ref{cor:coarse-trace-scaling-extension} and \eqref{formula:integral convergence},
\begin{align*}
 \|R_{\rho}-S_{\rho}\|_{L^{p}(\mathbb{R})}^{p}
  &\leq c_{\rho}^{p}\int_{\mathbb{R}}\|\mathcal{F}_{G,p}\Big(f_{(\cdot)} e^{-\rho((\cdot)-s)^{2}}-f_{s} e^{-\rho((\cdot)-s)^{2}}\Big)\|_{\mathcal{F}L^{q}(G)}^{p}ds\\
&\leq c_{\rho}^{p}\int_{\mathbb{R}}\int_{\mathbb{R}}e^{-\omega t}e^{-p\rho(t-s)^{2}}\|f_{t}-f_{s}\|_{L^{p}(N)}^{p}\,dt\,ds\\
  &=\left(\frac{p\rho}{\pi}\right)^{\frac{1}{2}}\int_{\mathbb{R}}e^{-p\rho r^{2}}D_{f}(r)\,dr
  \longrightarrow0\quad\mbox{as}\quad \rho\rightarrow\infty.
\end{align*}
Combining this estimate with equality $\mathcal{L}_{\rho}(f)^{\frac{1}{p}}=\|R_{\rho}\|_{L^{p}(\mathbb{R})}$ and \eqref{eq:unified-frozen-field-bound},
\begin{align*}
\limsup_{\rho\rightarrow\infty}\mathcal{L}_{\rho}(f)^{\frac{1}{p}}
=\limsup_{\rho\rightarrow\infty}\|R_{\rho}\|_{L^{p}(\mathbb{R})}
=\limsup_{\rho\rightarrow\infty}\|S_{\rho}\|_{L^{p}(\mathbb{R})}
\leq B_{p}H_{p}(N)\|f\|_{L^{p}(G)}.
\end{align*}
The proof is complete.
\end{proof}

To show Theorem \ref{thm:one-dimensional-extension}, we firstly state the following result. 
\begin{lemma}\label{lem:H1-H2-endpoints}
For every second countable locally compact group $G$, one has
\begin{align*}
 H_{1}(G)=H_{2}(G)=1.
\end{align*}
\end{lemma}
The assertion for $p=2$ is Plancherel's identity.  The assertion for
$p=1$ follows from the contractive integrated representation and a
positive normalized $L^{1}$-approximate identity; see
\cite[Chapter~2]{Folland2016}.

Let us prove Theorem \ref{thm:one-dimensional-extension}.
\begin{proof}[Proof of Theorem~\ref{thm:one-dimensional-extension}]
Let $G=N\rtimes_{\alpha}\mathbb{R}$ and let $f\in L^{p}(G)$. For $1<p<2$, it follows from Lemmas \ref{lem:group-crossed-product-identification}, \ref{lem:flow-window-density} and \ref{lem:trace-preserving-freezing} that
\begin{align*}
\|\mathcal{F}_{G,p}f\|_{\mathcal{F}L^{q}(G)}
  =\lim_{\rho\downarrow0}\mathcal{L}_{\rho}(f)^{\frac{1}{p}}
  \leq\limsup_{\rho\to\infty}\mathcal{L}_{\rho}(f)^{\frac{1}{p}}
  \leq B_{p}H_{p}(N)\|f\|_{L^{p}(G)}.
\end{align*}
At $p=1$ and $p=2$, Lemma~\ref{lem:H1-H2-endpoints} gives
\begin{align*}
 H_{1}(G)=H_{1}(N)=B_{1}=1,
 \qquad
 H_{2}(G)=H_{2}(N)=B_{2}=1.
\end{align*}
Thus the endpoint cases follow as well.  The proof is complete.
\end{proof}

\begin{remark}
Every formula in this section is valid for $\omega=0$.  In that case all exponential weights are one, the dual weight is a trace, and the exponential factor is one; hence the proof reduces to the trace-preserving crossed-product flow.  
\end{remark}

\section{Completely solvable Lie groups}\label{sec:completely-solvable-proof}
In this section, we prove Theorem \ref{thm:completely-solvable-sharp}. Before that, we should investigate some properties for solvable Lie groups.


The standard ideal-flag structure needed here is collected in
Lemma~\ref{lem:completely-solvable-ideal-flag}.  The fact that every
member of the flag is an ideal of the final algebra is essential in the
following trace argument.
\begin{lemma}\label{lem:unimodularity-along-ideal-flag}
Let $\mathfrak{g}$ be a unimodular Lie algebra and let $\mathfrak{h}\subset\mathfrak{g}$ be an ideal.  Then
$\mathfrak{h}$ is unimodular.  Consequently, if Lemma \ref{lem:completely-solvable-ideal-flag} \eqref{lem:flag decomposition} is a complete ideal flag in a unimodular completely solvable algebra, then every $\mathfrak{g}_{j}$ is
unimodular, and the action $\alpha_{j}$ in Lemma \ref{lem:completely-solvable-ideal-flag} \eqref{lem:flag isomorphism} preserves Haar measure on $G_{j-1}$.
\end{lemma}
\begin{proof}
For $X\in\mathfrak{h}$, the map induced by $\operatorname{ad}_{\mathfrak{g}}X$ on $\mathfrak{g}/\mathfrak{h}$ is zero, because $[X,\mathfrak{g}]\subseteq\mathfrak{h}$.  Therefore
\begin{align*}
 \operatorname{tr}\bigl(\operatorname{ad}_{\mathfrak{h}}X\bigr)
 = \operatorname{tr}\bigl(\operatorname{ad}_{\mathfrak{g}}X\bigr)=0.
\end{align*}
Thus $\mathfrak{h}$ is unimodular.   Applying it with $\mathfrak{h}=\mathfrak{g}_{j}$, if
$X_{j}\in\mathfrak{g}_{j}\setminus\mathfrak{g}_{j-1}$, then $\operatorname{ad}X_{j}$ acts trivially on the one-dimensional quotient $\mathfrak{g}_{j}/\mathfrak{g}_{j-1}$, whence
$\operatorname{tr}\big(\operatorname{ad}X_{j}\big|_{\mathfrak{g}_{j-1}}\big)
 = \operatorname{tr}\big(\operatorname{ad}X_{j}\big|_{\mathfrak{g}_{j}}\big)=0.$
It follows that
\begin{align*}
 \big|\det \big(d\alpha_{j}(t)\big|_{\mathfrak{g}_{j-1}}\big)\big|
 =\exp\big(t\operatorname{tr}(\operatorname{ad}X_{j}\big|_{\mathfrak{g}_{j-1}}\big)\big)
 =1.
\end{align*}
Therefore, $\alpha_{j}$ preserves Haar measure on $G_{j-1}$.
\end{proof}

\begin{lemma}\label{lem:unimodular-modular-kernel}
Suppose that $G$ is simply connected, nonunimodular and completely solvable.   
\begin{enumerate}[\rm(i)]
\item\label{lem:assertion ideal} Then $\mathfrak{n}:=\ker\delta_{\mathfrak{g}}$ is a completely solvable, unimodular ideal of codimension one. 

\item\label{lem:kernel of N} The analytic subgroup $N$ with Lie algebra $\mathfrak{n}$ is $\ker\Delta_{G}$; it is closed, connected, simply connected, and unimodular. 

\item\label{lem:cross product decomposition} If $X_{0}\in\mathfrak{g}$ is chosen so that $\omega=\delta_{\mathfrak{g}}(X_{0})>0,$ then
\begin{align*}
 G\simeq N\rtimes_{\alpha}\mathbb{R},
 \qquad\alpha_{t}(x)=\exp(tX_{0})x\exp(-tX_{0}),t\in\mathbb{R},x\in N,
\end{align*}
and, for a Haar measure $dx$ on $N$,
\begin{align*}
 d(\alpha_{t} x)=e^{\omega t}\,dx,\qquad d\mu_{G}(x,t)=e^{-\omega t}\,dx\,dt,
 \qquad \Delta_{G}(x,t)=e^{-\omega t}.
\end{align*}
\end{enumerate} 
\end{lemma}
\begin{proof}
Clearly,
\begin{align*}
 \delta_{\mathfrak{g}}([X,Y])=\operatorname{tr}([\operatorname{ad}_{\mathfrak{g}}X,\operatorname{ad}_{\mathfrak{g}}Y])=0.
\end{align*}
Thus $\delta_{\mathfrak{g}}$ is a Lie algebra character.  Since it is nonzero, ${\rm ker}(\delta_{\mathfrak{g}})$ is a completely solvable ideal of codimension one.  If $X\in\mathfrak{n}$, then $\operatorname{ad}_{\mathfrak{g}}X$ induces the zero map on $\mathfrak{g}/\mathfrak{n}$; consequently, $\operatorname{tr}\bigl(\operatorname{ad}_{\mathfrak{n}}X\bigr)=\operatorname{tr}\bigl(\operatorname{ad}_{\mathfrak{g}}X\bigr)=0.$ Hence $\mathfrak{n}$ is unimodular and \eqref{lem:assertion ideal} is claimed.

Identifying $N$ with $\exp(\mathfrak{n})$, one has
\begin{align*}
 \Delta_{G}(\exp X)=\exp\bigl(-\delta_{\mathfrak{g}}(X)\bigr),\qquad X\in\mathfrak{g}.
\end{align*}
By the exponential and codimension-one splitting facts recalled after
Lemma~\ref{lem:completely-solvable-ideal-flag}, it follows that
$N=\ker\Delta_{G}$, that $N$ is closed, connected, and simply
connected, and $G\simeq N\rtimes_{\alpha}\mathbb{R}$ with $\alpha_{t}$ given in \eqref{lem:cross product decomposition}.  Since
$\operatorname{ad}X_{0}$ acts trivially on
$\mathfrak{g}/\mathfrak{n}$,
\begin{align*}
 \operatorname{tr}
 \left(\left.\operatorname{ad}X_{0}\right|_{\mathfrak{n}}\right)
 =\delta_{\mathfrak{g}}(X_{0})=\omega.
\end{align*}
Taking determinants of $d\alpha_{t}$ gives the measure and modular function in \eqref{lem:cross product decomposition}.
\end{proof}


\begin{proof}[Proof of Theorem \ref{thm:completely-solvable-sharp}]
Assume first that $1<p<2$ and let $q=\frac{p}{p-1}$.  Suppose that
$G$ is unimodular.  Choose the complete ideal flag
$\{\mathfrak{g}_{j}\}_{j=0}^{\dim{G}}$ from Lemma
\ref{lem:completely-solvable-ideal-flag}.  By Lemma
\ref{lem:unimodularity-along-ideal-flag}, every associated group $G_{j}$
is unimodular and every action $\alpha_{j}$ preserves the corresponding
Haar measure.

Let $\tau_{j}$ denote the Plancherel trace of $G_{j}$.  At the step
\begin{align*}
 G_{j}=G_{j-1}\rtimes_{\alpha_{j}}\mathbb{R},
\end{align*}
the dual weight $\Phi_{j}$ is the Plancherel trace $\tau_{j}$ because
$\alpha_{j}$ preserves Haar measure.  Using notation in Lemma \ref{lem:weighted-crossed-product-assembly}, if $f$ belongs to the
corresponding Fourier--Tomita core and $f_{t}(x)=f(x,t)$,  then
\begin{align*}
 \mathcal{A}_{0}\left(
 t\longmapsto
 j_{q}^{\tau_{j-1}}\bigl(\lambda_{G_{j-1}}(f_{t})\bigr)
 \right)
 =j_{q}^{\Phi_{j}}\bigl(\lambda_{G_{j}}(f)\bigr)
 =j_{q}^{\tau_{j}}\bigl(\lambda_{G_{j}}(f)\bigr).
\end{align*}
Equivalently,
\begin{align*}
 \mathcal{A}_{0}\left(
 t\longmapsto\mathcal{F}_{G_{j-1},p}f_{t}
 \right)
 =\mathcal{F}_{G_{j},p}f
 \quad\text{in}\quad\mathcal{F}L^{q}(G_{j}).
\end{align*}
Applying Theorem \ref{thm:one-dimensional-extension} with $\omega=0$, expressed through these
identifications, one has
\begin{align*}
 \|\mathcal{F}_{G_{j},p}f\|_{\mathcal{F}L^{q}(G_{j})}
 \leq B_{p}H_{p}(G_{j-1})\|f\|_{L^{p}(G_{j})}.
\end{align*}
Starting with the trivial group $G_{0}$, for which $H_{p}(G_{0})=1$,
and taking the supremum over nonzero $f$ at every step, it follows that
\begin{align*}
 H_{p}(G_{j})\leq B_{p}H_{p}(G_{j-1})\leq B_{p}^{j},
 \qquad 1\leq j\leq\dim{G}.
\end{align*}
In particular, $H_{p}(G)\leq B_{p}^{\dim{G}}$.  Together with the
standard local lower bound \eqref{formula:lower bound}, this proves that,
in the unimodular setting, $H_{p}(G)=B_{p}^{\dim{G}}$.

When $G$ is nonunimodular, using Lemma \ref{lem:unimodular-modular-kernel},
\begin{align*}
 G=N\rtimes_{\alpha}\mathbb{R},\qquad
 \left|\det((d\alpha_{t})_{e})\right|=e^{\omega t},\qquad \omega>0,
\end{align*}
where $N$ is simply connected, unimodular, and completely solvable with
$\dim{N}=\dim{G}-1$.  Let $\tau_{N}$ be its Plancherel trace, let
$\beta$ be the induced action on $\operatorname{VN}(N)$, and let
$\Phi_{G}$ be the dual weight on
$\operatorname{VN}(G)
 \simeq\operatorname{VN}(N)\rtimes_{\beta}\mathbb{R}.$
Then $\Phi_{G}$ is the
Plancherel weight of $G$.  On the Fourier--Tomita core,
\begin{align*}
 j_{q}^{\Phi_{G}}\bigl(\lambda_{G}(f)\bigr)
 =j_{q}^{\Phi_{G}}\left(
 \int_{\mathbb{R}}e^{-\omega t}
 \pi_{\beta}\bigl(\lambda_{N}(f_{t})\bigr)\lambda_{t}\,dt
 \right)
 =\mathcal{A}_{\omega}\left(
 t\longmapsto
 j_{q}^{\tau_{N}}\bigl(\lambda_{N}(f_{t})\bigr)
 \right).
\end{align*}
Consequently,
\begin{align*}
 \|\mathcal{F}_{G,p}f\|_{\mathcal{F}L^{q}(G)}
 =\|j_{q}^{\Phi_{G}}(\lambda_{G}(f))\|_{L^{q}(\operatorname{VN}(G))_{\mathrm{left}}}\quad\mbox{and}\quad
 \|\mathcal{F}_{N,p}f_{t}\|_{\mathcal{F}L^{q}(N)}
 =\|j_{q}^{\tau_{N}}(\lambda_{N}(f_{t}))\|_{L^{q}(\operatorname{VN}(N),\tau_{N})}.
\end{align*}
The index is $q$, not $p$, because the concrete Fourier-side space is
obtained at interpolation parameter $\frac{2}{q}$.  Applying the result
already proved for the unimodular group $N$ gives
\begin{align*}
 H_{p}(N)=B_{p}^{\dim{G}-1}.
\end{align*}
Applying Theorem \ref{thm:one-dimensional-extension} with $\omega>0$,
\begin{align*}
 \|\mathcal{F}_{G,p}f\|_{\mathcal{F}L^{q}(G)}
 \leq B_{p}H_{p}(N)\|f\|_{L^{p}(G)},
\end{align*}
and hence
\begin{align*}
H_{p}(G)\leq B_{p}H_{p}(N)=B_{p}^{\dim{G}}.
\end{align*}
Combining the local lower bound \eqref{formula:lower bound}, the assertion for $1<p<2$ follows.

The endpoint cases follow from Lemma~\ref{lem:H1-H2-endpoints}:
$H_{1}(G)=H_{2}(G)=1$, while
$B_{1}^{\dim G}=B_{2}^{\dim G}=1$.  The proof is complete.
\end{proof}

\begin{remark}\label{rem:modular-kernel-order}
An arbitrary codimension-one flag in a nonunimodular algebra need not make the intermediate actions trace preserving.  The filtration used above is deliberately different: one first takes a complete ideal flag inside $\ker\delta_{\mathfrak{g}}$, so that every lower-dimensional step is trace preserving, and only then adjoins the single modular direction. 
\end{remark}

\section{The finite-center semisimple theorem}\label{sec:semisimple-proof}

In this section we prove Theorem~\ref{thm:semisimple-sharp}.  
We employ the Cartan decomposition and the Iwasawa decomposition for semisimple Lie groups. To use these tricks, we first need to establish a closed-subgroup comparison principle for sharp constants.

\subsection{Closed-subgroup comparison principle}
\label{subsec:closed-subgroup-comparison}

We use throughout the Fourier--Lebesgue, Tomita-core, and left-coset
conventions collected in
Section~\ref{subsec:preliminary-group-Fourier-transform} and
Lemma~\ref{lem:left-coset-disintegration}.

The next lemma is the point at which nonunimodularity of $H$ must be treated
carefully.  We use the inclusion of group von Neumann algebras only at the
operator-norm endpoint.  We do not claim that this inclusion preserves the
Plancherel weights.

\begin{lemma}\label{lem:compatible-endpoint-assembly}
Let $G$ be a second countable locally compact group and let $H$ be a
closed subgroup. There are compatible contractions
\begin{align*}
\mathcal{A}_{\infty}:L^{1}(H_{G};\operatorname{Cv}^{2}(H))\longrightarrow \operatorname{Cv}^{2}(G),\qquad
\mathcal{A}_{2}:L^{2}(H_{G};L^{2}(H))\longrightarrow L^{2}(G).
\end{align*}
If $0<\theta<1$ and $\frac{1}{u}=1-\frac{\theta}{2}$, the complex interpolation of $\mathcal{A}_{\infty}$ and $\mathcal{A}_{2}$ is a contraction
\begin{align*}
\mathcal{A}_{u,\frac{2}{\theta}}:
L^{u}(H_{G};\mathcal{F}L^{\frac{2}{\theta}}(H))
\longrightarrow\mathcal{F}L^{\frac{2}{\theta}}(G).
\end{align*}
\end{lemma}
\begin{proof}
Define
\begin{align*}
    \iota_{H}:\operatorname{VN}(H)\longrightarrow\operatorname{VN}(G),\qquad
    \iota_{H}(\lambda_{H}(h))=\lambda_{G}(h).
\end{align*}
By Lemma \ref{lem:left-coset-disintegration} \eqref{unitary adjoint represation},  $\iota_{H}$ is a normal isometric embedding. For $u\in L^{1}(H_{G};\operatorname{Cv}^{2}(H))$, define the ultraweak integral
\begin{align}
\mathcal{A}_{\infty} u=\int_{H_{G}}^{\mathrm{uw}}\iota_{H}(u(x))\lambda_{G}(s(x))\,d\nu(x).
    \label{eq:A-infinity-definition}
\end{align}
Since right multiplication by a unitary preserves the operator norm, it follows that
\begin{align*}
\|\mathcal{A}_{\infty} u\|_{\operatorname{Cv}^{2}(G)}\leq\int_{H_{G}}\|u(x)\|_{\operatorname{Cv}^{2}(H)}\,d\nu(x).
\end{align*}
At the Hilbert-space endpoint, define
\begin{align*}
    \Big(\mathcal{A}_{2}u\Big)(hs(x))=u(x)(h),\quad u\in L^{2}(H_{G};L^{2}(H)),\;x\in H_{G},\;h\in H.
\end{align*}
By Lemma \ref{lem:left-coset-disintegration} \eqref{coset representation}, $\mathcal{A}_{2}$ is an isometry.

We verify that these two maps represent the same distribution.  Take
$\phi(x)=\sum_{j=1}^{m}\alpha_{j}(x)a_{j}$ with
$a_{j}\in C_{c}^{\infty}(H)$, where each $\alpha_{j}$ is a bounded
simple function of finite-measure support.  At the operator endpoint,
$a_{j}$ is identified with $\lambda_{H}(a_{j})$.  Let
\begin{align*}
F_{\phi}(hs(x))=\sum_{j=1}^{m}\alpha_{j}(x)a_{j}(h).
\end{align*}
Clearly, $\mathcal{A}_{2}\phi=F_{\phi}$ and
\begin{align*}
\mathcal{A}_{\infty}\phi=\int_{H_{G}}\int_{H} F_{\phi}(hs(x))\lambda_{G}(hs(x))\,dh\,d\nu(x)=\lambda_{G}(F_{\phi}).
\end{align*}
Then the compatibility extends to all of $L^{1}(H_{G};\operatorname{Cv}^{2}(H))\cap L^{2}(H_{G};L^{2}(H))$, not only to the algebraic core.  If $k^{\phi}_{x}\in L^{2}(H)$ is the convolution kernel
of $\phi(x)$, set $F_{\phi}(hs(x))=k^{\phi}_{x}(h)$.  For $\xi,\eta\in C_{c}^{\infty}(G)$, the coefficient $h\longmapsto \langle\lambda_{G}(hs(x))\xi,\eta\rangle$ belongs to the Fourier algebra $A(H)$.  Pairing
\eqref{eq:A-infinity-definition} with this coefficient,
\begin{align*}
    \langle\mathcal{A}_{\infty} \phi\,\xi,\eta\rangle
    =\int_{H_{G}}\int_{H} k^{\phi}_{x}(h)\langle\lambda_{G}(hs(x))\xi,\eta\rangle \,dh\,d\nu(x)
    =\int_{G} F_{\phi}(g)\langle\lambda_{G}(g)\xi,\eta\rangle\,dg.
\end{align*}
Thus, $\mathcal{A}_{\infty}$ and $\mathcal{A}_{2}$ agree in the ambient
distribution space.  Let $\mathcal{A}_{u,\frac{2}{\theta}}$ be the
interpolated extension, where
$\frac{1}{u}=1-\frac{\theta}{2}$.  For $\frac{1}{u}=1-\frac{\theta}{2}$, vector-valued complex interpolation
\cite{BerghLofstrom1976,PisierXu2003} asserts that
\begin{align*}
 \left[
 L^{1}(H_{G};\operatorname{Cv}^{2}(H)),
 L^{2}(H_{G};L^{2}(H))
 \right]_{\theta}
 &=
 L^{u}\left(
 H_{G};
 [\operatorname{Cv}^{2}(H),L^{2}(H)]_{\theta}
 \right)\\
 &=
 L^{u}\left(H_{G};\mathcal{F}L^{\frac{2}{\theta}}(H)\right).
\end{align*}
Here Lemma~\ref{lem:semisimple-Tomita-core} supplies the common dense
compatible core, and \eqref{FL interpolation 2-infty} identifies the
inner interpolation space and the target space.  Thus the source and
target are precisely the spaces stated in the lemma.

We finally record its concrete form.  At the two endpoints, the
compatibility proved above says, under the canonical endpoint
identifications,
\begin{align*}
 \mathcal{A}_{\infty}\left(
 x\longmapsto
 j_{\infty}^{\varphi_{H}}(\lambda_{H}(\phi(x)))
 \right)
 =j_{\infty}^{\varphi_{G}}(\lambda_{G}(F_{\phi})),\quad
 \mathcal{A}_{2}\left(
 x\longmapsto
 j_{2}^{\varphi_{H}}(\lambda_{H}(\phi(x)))
 \right)
 =j_{2}^{\varphi_{G}}(\lambda_{G}(F_{\phi})).
\end{align*}
The second equality is understood through the two Plancherel GNS
identifications.  If $1<p<2$ and $q=\frac{p}{p-1}$, interpolation and
transport through the canonical Fourier--Lebesgue isometries give
\begin{align*}
 \mathcal{A}_{p,q}\left(
 x\longmapsto
 j_{q}^{\varphi_{H}}(\lambda_{H}(\phi(x)))
 \right)
 =j_{q}^{\varphi_{G}}(\lambda_{G}(F_{\phi})).
\end{align*}
Here the same notation denotes the transported map
\begin{align*}
 \mathcal{A}_{p,q}:
 L^{p}\bigl(H_{G};L^{q}(\operatorname{VN}(H))_{\mathrm{left}}\bigr)
 \longrightarrow L^{q}(\operatorname{VN}(G))_{\mathrm{left}}.
\end{align*}
This equality follows from compatible interpolation of the two endpoint
maps; it does not assert that $\iota_{H}$ preserves the Plancherel
weights or directly induces an $L^{q}$ embedding.
\end{proof}

\begin{theorem}\label{thm:closed-subgroup-comparison}
Let $G$ be a second countable Lie group and let $H$ be a closed subgroup.  We have
\begin{align*}
H_{p}(G)\leq H_{p}(H), \qquad 1\leq p\leq2.
\end{align*}
\end{theorem}
\begin{proof}
Suppose first that $1<p<2$ and let $q=\frac{p}{p-1}$.  For
$f\in C_{c}^{\infty}(G)$, define
\begin{align*}
    f_{x}(h)=f(hs(x)),\quad x\in H_{G},\;h\in H.
\end{align*}
By the Hausdorff--Young inequality on $H$ and Lemma \ref{lem:left-coset-disintegration} \eqref{coset representation},
\begin{align*}
 \int_{H_{G}}
 \|\mathcal{F}_{H,p}f_{x}\|_{\mathcal{F}L^{q}(H)}^{p}\,d\nu(x)
 \leq H_{p}(H)^{p}
 \int_{H_{G}}\|f_{x}\|_{L^{p}(H)}^{p}\,d\nu(x)
 =H_{p}(H)^{p}\|f\|_{L^{p}(G)}^{p}.
\end{align*}
Lemma \ref{lem:compatible-endpoint-assembly} gives the abstract
Fourier--Lebesgue identity
\begin{align*}
 \mathcal{A}_{p,q}\left(
 x\longmapsto\mathcal{F}_{H,p}f_{x}
 \right)
 =\mathcal{F}_{G,p}f
 \quad\text{in}\quad\mathcal{F}L^{q}(G).
\end{align*}
On the Fourier--Tomita core, its concrete realization is
\begin{align*}
 \mathcal{A}_{p,q}\left(
 x\longmapsto
 j_{q}^{\varphi_{H}}(\lambda_{H}(f_{x}))
 \right)
 =j_{q}^{\varphi_{G}}(\lambda_{G}(f)).
\end{align*}
Consequently,
\begin{align*}
 \|\mathcal{F}_{G,p}f\|_{\mathcal{F}L^{q}(G)}
 \leq\left(
 \int_{H_{G}}
 \|\mathcal{F}_{H,p}f_{x}\|_{\mathcal{F}L^{q}(H)}^{p}
 \,d\nu(x)
 \right)^{\frac{1}{p}}
 \leq H_{p}(H)\|f\|_{L^{p}(G)}.
\end{align*}
The result for $1<p<2$ follows by density.   At $p=1$
and $p=2$, Lemma~\ref{lem:H1-H2-endpoints} gives
$H_{1}(G)=H_{1}(H)=1$ and $H_{2}(G)=H_{2}(H)=1$.  The proof is complete.
\end{proof}

\subsection{Contraction to the Cartan motion group}\label{subsec:Cartan-contraction}

Under the assumptions of Theorem~\ref{thm:semisimple-sharp}, we use
the Iwasawa decomposition, Cartan decomposition and polar data fixed in
Section~\ref{subsec:preliminary-semisimple-Lie-groups}, including the measure
normalization and the Jacobian $J$ in \eqref{eq:Cartan-Jacobian}.
The contraction of a semisimple Lie group to its Cartan motion group is
classical; see \cite{DooleyRice1985}.  We use below the explicit polar
realization adapted to the present norm comparison.

For $0<\varepsilon\leq1$, set 
\begin{align}\label{K-p isomorphism}
\Phi_{\varepsilon}(k,X)=k\exp(\varepsilon X),\quad(k,X)\in K\times\mathfrak{p}.
\end{align}
Pulling the multiplication of $G$ back to the manifold $K\times\mathfrak{p}$, the resulting group is denoted by $G_{\varepsilon}=(K\times\mathfrak{p},*_{\varepsilon})$, equipped with multiplication
\begin{align*}
(k,X)*_{\varepsilon}(l,Y)
=\Phi_{\varepsilon}^{-1}\Big(\Phi_{\varepsilon}(k,X)\Phi_{\varepsilon}(l,Y)\Big).
\end{align*}
For each $\varepsilon>0$, $G_{\varepsilon}$ is isomorphic to $G$.
\begin{lemma}\label{lem:contracted-group-law-Haar}
Under the assumptions in Theorem~\ref{thm:semisimple-sharp}, let
$G_{0}=(K\times\mathfrak{p},*_{0})$ have multiplication, unit, and inversion
\begin{align*}
(k,X)*_{0}(l,Y)=\bigl(kl,\operatorname{Ad}(l^{-1})X+Y\bigr),\quad
e_{G_{0}}=(e_{K},0),\quad
(k,X)^{-1_{0}}=\bigl(k^{-1},-\operatorname{Ad}(k)X\bigr).
\end{align*}
We have the following results.
\begin{enumerate}[\rm(i)]
\item\label{G0 isomorphism} The map $(k,X)\longmapsto(k,\operatorname{Ad}(k)X)$ is a Haar-preserving isomorphism from $G_{0}$ onto the standard Cartan motion group $K\ltimes_{\operatorname{Ad}}\mathfrak{p}$ with multiplication $(k,X)(l,Y)=(kl,X+\operatorname{Ad}(k)Y).$

\item\label{G0 operations} As $\varepsilon\downarrow0$, multiplication and inversion on $G_{\varepsilon}$ converge smoothly on compact subsets to those on $G_{0}$.

\item\label{G0 measure convergence} Let
\begin{align*}
dm_{\varepsilon}
=\varepsilon^{-\dim\mathfrak{p}}\Phi_{\varepsilon}^{*}(dg)
=J(\varepsilon X)\,dk\,dX,
\end{align*}
where $\Phi_{\varepsilon}$ and $J$ are given in \eqref{K-p isomorphism} and \eqref{eq:Cartan-Jacobian}.
Then $dm_{\varepsilon}$ is the Haar measure on $G_{\varepsilon}$, and it converges in $C^{\infty}_{\mathrm{loc}}$ to $dm_{0}:=dk\,dX.$

\item\label{G0 unimodularity} For all $0<\varepsilon\leq1$, $G_{\varepsilon}$ and $G_{0}$ are unimodular.
\end{enumerate}
\end{lemma}
\begin{proof}
The asserted isomorphism in \eqref{G0 isomorphism} with the usual semidirect-product convention follows by direct calculations.

We have
\begin{align*}
k\exp(\varepsilon X)l\exp(\varepsilon Y)=kl\exp(\varepsilon\operatorname{Ad}(l^{-1})X)\exp(\varepsilon Y).
\end{align*}
Near $(0,0)\in\mathfrak{p}\times\mathfrak{p}$, write
\begin{align*}
\exp U\exp V=\kappa(U,V)\exp Z(U,V),\qquad \kappa(U,V)\in K,\quad Z(U,V)\in\mathfrak{p}.
\end{align*}
By the smooth polar decomposition and its differential at the origin,
\begin{align*}
\kappa(\varepsilon U,\varepsilon V)=e+O(\varepsilon^{2}),\qquad
Z(\varepsilon U,\varepsilon V)=\varepsilon(U+V)+O(\varepsilon^{2}),
\end{align*}
with all derivatives and remainders controlled on compact sets.  The inverse formula follows from
\begin{align*}
k^{-1}\exp(-\varepsilon\operatorname{Ad}(k)X)=\exp(-\varepsilon X)k^{-1}.
\end{align*}
This proves \eqref{G0 operations}.  

Observe that
\begin{align*}
\Phi_{\varepsilon} =\Theta\circ[(k,X)\mapsto(k,\varepsilon X)],
\end{align*}
The formula \eqref{eq:Cartan-Jacobian} gives \eqref{G0 measure convergence}. Since $G_{\varepsilon}\simeq G$ for $\varepsilon>0$, and $K$ is compact and acts orthogonally on $\mathfrak{p}$ for $\varepsilon=0$,  the unimodularity in \eqref{G0 unimodularity} follows.
\end{proof}

Denote by $\lambda_{\varepsilon}$ and $\tau_{\varepsilon}$ the left
regular representation and the Plancherel trace of $G_{\varepsilon}$,
respectively; use $\lambda_{0},\tau_{0}$ and $\lambda_{G},\tau_{G}$ for
$G_{0}$ and $G$.  All these groups are unimodular.  Therefore, for
$2\leq r\leq\infty$,
\begin{align*}
 \mathcal{F}L^{r}(G_{\varepsilon})
 \cong L^{r}(\operatorname{VN}(G_{\varepsilon}),\tau_{\varepsilon}),\quad
 \mathcal{F}L^{r}(G_{0})
 \cong L^{r}(\operatorname{VN}(G_{0}),\tau_{0}),\quad
 \mathcal{F}L^{r}(G)
 \cong L^{r}(\operatorname{VN}(G),\tau_{G}).
\end{align*}
For $1\leq p\leq2$ and $q=\frac{p}{p-1}$, with $q=\infty$ when $p=1$, the respective Fourier--Tomita core identifications are
\begin{align*}
 \mathcal{F}_{G_{\varepsilon},p}F
 \longleftrightarrow
 j_{q}^{\tau_{\varepsilon}}(\lambda_{\varepsilon}(F)),\quad
 \mathcal{F}_{G_{0},p}F
 \longleftrightarrow j_{q}^{\tau_{0}}(\lambda_{0}(F)),\quad
 \mathcal{F}_{G,p}f
 \longleftrightarrow j_{q}^{\tau_{G}}(\lambda_{G}(f)).
\end{align*}
For $u_{\varepsilon},v_{\varepsilon}\in C_{c}(K\times\mathfrak{p})$,
define
\begin{align}\label{G-varepsilon convolution}
(u_{\varepsilon} *_{\varepsilon} v_{\varepsilon})(x)=
\int_{K\times\mathfrak{p}}u_{\varepsilon}(y)v_{\varepsilon}(y^{-1_{\varepsilon}} *_{\varepsilon} x)\,dm_{\varepsilon}(y),\quad x\in K\times\mathfrak{p}.
\end{align}
The next two lemmas extend the elementary convergence of integer convolution moments to every real exponent $q\geq2$.
\begin{lemma}\label{lem:contraction-integer-moments}
Under the assumptions in Theorem~\ref{thm:semisimple-sharp}, let
$G_{\varepsilon}$ be the contracted groups defined above.  Let
$u_{\varepsilon},v_{\varepsilon}\in C_{c}(K\times\mathfrak{p})$ have supports in fixed compact
sets and converge uniformly to $u_{0},v_{0}\in C_{c}(K\times\mathfrak{p})$.  Then
\begin{align*}
    u_{\varepsilon} *_{\varepsilon} v_{\varepsilon}
    \longrightarrow u_{0} *_{0} v_{0},\quad\mbox{as}\quad\varepsilon\rightarrow0,
\end{align*}
uniformly, and all convolutions have support in one compact set.  The same holds for every fixed iterated convolution and for convolution involution.

In particular, if $F\in C_{c}^{\infty}(K\times\mathfrak{p})$ and $A_{F,\varepsilon}:=(\lambda_{\varepsilon}(F))^{*}\lambda_{\varepsilon}(F)$, then there is $c_{F}<\infty$, independent of $0<\varepsilon\leq1$, such that
$0\leq A_{F,\varepsilon}\leq c_{F}{\rm id}$, $\tau_{\varepsilon}(A_{F,\varepsilon})<\infty,$ and
\begin{align*}
\tau_{\varepsilon}(A_{F,\varepsilon}^{n})\longrightarrow \tau_{0}(A_{F,0}^{n}),
\quad\mbox{as}\quad\varepsilon\rightarrow0,
\end{align*}
for every integer $n\geq1$.
\end{lemma}
\begin{proof}
By Lemma \ref{lem:contracted-group-law-Haar} \eqref{G0 measure convergence}, we have the uniform convergence of the integrand in \eqref{G-varepsilon convolution} on the relevant compact set.  The union, over $0\leq \varepsilon\leq1$, of the products of the two fixed support sets is compact. The assertion $u_{\varepsilon} *_{\varepsilon} v_{\varepsilon}\rightarrow u_{0} *_{0} v_{0}$ and its iterated versions follow from the dominated convergence theorem.

Note that  $A_{F,\varepsilon}=\lambda_{\varepsilon}(F^{*_{\varepsilon}} *_{\varepsilon}F)$. Using Young's inequality,
\begin{align*}
\|A_{F,\varepsilon}\|_{\infty}=\|\lambda_{\varepsilon}(F^{*_{\varepsilon}} *_{\varepsilon}F)\|_{\infty}
\leq\sup_{0<\varepsilon\leq1}\|F\|_{L^{1}(K\times\mathfrak{p},dm_{\varepsilon})}^{2}=:c_{F}.
\end{align*}
By the Plancherel identity and the unimodularity of $G_{\varepsilon}$ in Lemma \ref{lem:contracted-group-law-Haar}  \eqref{G0 unimodularity}, one has
\begin{align*}
\tau_{\varepsilon}(A_{F,\varepsilon}^{n}) =(F^{*_{\varepsilon}} *_{\varepsilon}F)^{(*_{\varepsilon} n)}(e),
    \qquad n\geq1.
\end{align*}
Then the convergence for $\tau_{\varepsilon}(A_{F,\varepsilon}^{n})$ follows from the established convergence of the iterated convolutions.
\end{proof}

\begin{lemma}\label{lem:weighted-spectral-measures}
Under the assumptions in Theorem~\ref{thm:semisimple-sharp}, let
$G_{\varepsilon}$ be the contracted groups defined above and let
$F\in C_{c}^{\infty}(K\times\mathfrak{p})$. For every $2\leq q<\infty$, let $p=\frac{q}{q-1}$. Then
\begin{align*}
\|\mathcal{F}_{G_{\varepsilon},p}F\|_{\mathcal{F}L^{q}(G_{\varepsilon})}
\longrightarrow
\|\mathcal{F}_{G_{0},p}F\|_{\mathcal{F}L^{q}(G_{0})}
\quad\mbox{as}\quad\varepsilon\rightarrow0.
\end{align*}
\end{lemma}

\begin{proof}
Let $A_{F,\varepsilon}$ and the constant $c_{F}$ be given in Lemma \ref{lem:contraction-integer-moments}. Let $E_{\varepsilon}$ be the spectral measure of $A_{F,\varepsilon}$.  The ordinary spectral distribution has infinite mass because $\tau_{\varepsilon}(I)$ is generally infinite.  Define a finite positive measure on $[0,c_{F}]$ by
\begin{align*}
\nu_{\varepsilon}(B)=\tau_{\varepsilon}(A_{F,\varepsilon} E_{\varepsilon}(B))\quad\mbox{for Borel set}\quad B\subset[0,c_{F}].
\end{align*}
Its moments satisfy
\begin{align*}
\int_{0}^{c_{F}}t^{k}\,d\nu_{\varepsilon}(t)=\tau_{\varepsilon}(A_{F,\varepsilon}^{k+1}),\qquad k=0,1,2,\ldots.
\end{align*}
By Lemma \ref{lem:contraction-integer-moments}, the convergence of all moments, including the total masses, follows.  Since polynomials are dense on continuous functions of $[0,c_{F}]$, it gives $\nu_{\varepsilon}\rightarrow\nu_{0}$ in the sense of weak topology. Since $t\mapsto t^{\frac{q}{2}-1}$ is continuous on $[0,c_{F}]$ for $q\geq2$, one has
\begin{align*}
\|\mathcal{F}_{G_{\varepsilon},p}F\|_{\mathcal{F}L^{q}(G_{\varepsilon})}^{q}
&=\|j_{q}^{\tau_{\varepsilon}}(\lambda_{\varepsilon}(F))\|_{L^{q}(\operatorname{VN}(G_{\varepsilon}),\tau_{\varepsilon})}^{q}\\
&=\tau_{\varepsilon}(A_{F,\varepsilon}^{\frac{q}{2}})\\
&=\int_{0}^{c_{F}}t^{\frac{q}{2}-1}\,d\nu_{\varepsilon}(t)
\longrightarrow\int_{0}^{c_{F}}t^{\frac{q}{2}-1}\,d\nu_{0}(t)\quad(\mbox{as}\quad\varepsilon\rightarrow0)\\
&=\|j_{q}^{\tau_{0}}(\lambda_{0}(F))\|_{L^{q}(\operatorname{VN}(G_{0}),\tau_{0})}^{q}\\
&=\|\mathcal{F}_{G_{0},p}F\|_{\mathcal{F}L^{q}(G_{0})}^{q}.
\end{align*}
The proof is complete.
\end{proof}

We next identify exactly how the Plancherel trace changes under the isomorphism $\Phi_{\varepsilon}$.

\begin{lemma}\label{lem:exact-Plancherel-scaling}
Under the assumptions in Theorem~\ref{thm:semisimple-sharp}, let
$\varepsilon>0$ and let $\Phi_{\varepsilon}$ be given in \eqref{K-p isomorphism}.  Let $1<p\leq2$ and $q=\frac{p}{p-1}$. If $F\in L^{p}(G_{\varepsilon})$ and $f_{F,\varepsilon}(\Phi_{\varepsilon}(k,X)):=\varepsilon^{-\frac{\dim(AN)}{p}}F(k,X)$, then
\begin{align*}
\|f_{F,\varepsilon}\|_{L^{p}(G)}
=\|F\|_{L^{p}(G_{\varepsilon},dm_{\varepsilon})},\quad
\|\mathcal{F}_{G,p}f_{F,\varepsilon}\|_{\mathcal{F}L^{q}(G)}
=\|\mathcal{F}_{G_{\varepsilon},p}F\|_{\mathcal{F}L^{q}(G_{\varepsilon})}.
\end{align*}
\end{lemma}
\begin{proof}
Let $c_{\varepsilon}=\varepsilon^{-\dim(AN)}$ and let
$F\in C_{c}^{\infty}(G_{\varepsilon})$.  Lemma
\ref{lem:contracted-group-law-Haar} \eqref{G0 measure convergence} and
the equality $\dim(AN)=\dim\mathfrak{p}$ from
Lemma~\ref{lem:AN-completely-solvable} give
\begin{align}\label{formula:Phi-dg}
(\Phi_{\varepsilon})_{*}(dm_{\varepsilon})=c_{\varepsilon}\,dg.
\end{align}
It directly gives that $\|f_{F,\varepsilon}\|_{L^{p}(G)}=\|F\|_{L^{p}(G_{\varepsilon},m_{\varepsilon})}$.

Define the operator $W_{\varepsilon}:L^{2}(G_{\varepsilon},m_{\varepsilon})\longrightarrow L^{2}(G,dg)$ by setting
$(W_{\varepsilon}\xi)(g)=c_{\varepsilon}^{\frac{1}{2}}\xi(\Phi_{\varepsilon}^{-1}g)$. By \eqref{formula:Phi-dg}, the operator $W_{\varepsilon}$ is unitary.
Set $\widetilde{F}_{\varepsilon}=F\circ\Phi_{\varepsilon}^{-1}$. Conjugating the regular representations gives that
\begin{align*}
W_{\varepsilon}\lambda_{\varepsilon}(F)W_{\varepsilon}^{-1}=c_{\varepsilon}\lambda_{G}(\widetilde{F}_{\varepsilon}).
\end{align*}
Combining this equality with Plancherel's theorem and \eqref{formula:Phi-dg}, one has
\begin{align*}
\tau_{\varepsilon}(\lambda_{\varepsilon}(F)^{*}\lambda_{\varepsilon}(F))
=\|F\|_{L^{2}(m_{\varepsilon})}^{2}
\end{align*}
and
\begin{align*}
\tau_{G}(W_{\varepsilon}\lambda_{\varepsilon}(F)^{*}\lambda_{\varepsilon}(F)W_{\varepsilon}^{-1})
=c_{\varepsilon}^{2}\|\widetilde{F}_{\varepsilon}\|_{L^{2}(G)}^{2}
=c_{\varepsilon}\|F\|_{L^{2}(m_{\varepsilon})}^{2}.
\end{align*}
Since $\lambda_{\varepsilon}(C_{c}^{\infty}(G_{\varepsilon}))$ is a left-Hilbert-algebra core for the Plancherel trace, the two normal faithful semifinite traces agree on the core squares. Uniqueness therefore gives the identity on $\operatorname{VN}(G_{\varepsilon})_{+}$, and spectral truncation extends it to every positive measurable operator $S$ affiliated with $\operatorname{VN}(G_{\varepsilon})$. Thus,
\begin{align}\label{eq:Plancherel-trace-scaling}
\tau_{\varepsilon}(S)=c_{\varepsilon}^{-1}\tau_{G}(W_{\varepsilon} S W_{\varepsilon}^{-1}).
\end{align}
Observe that
\begin{align*}
    \lambda_{G}(f_{F,\varepsilon})
    =c_{\varepsilon}^{\frac{1}{p}}\lambda_{G}(\widetilde{F}_{\varepsilon})
    =c_{\varepsilon}^{-\frac{1}{q}} W_{\varepsilon}\lambda_{\varepsilon}(F)W_{\varepsilon}^{-1}.
\end{align*}
Under the concrete tracial realizations of the Fourier--Lebesgue spaces, this gives
\begin{align*}
j_{q}^{\tau_{G}}(\lambda_{G}(f_{F,\varepsilon}))
=c_{\varepsilon}^{-\frac{1}{q}}
j_{q}^{\tau_{G}}
\left(W_{\varepsilon}\lambda_{\varepsilon}(F)W_{\varepsilon}^{-1}\right).
\end{align*}
Taking the $q$-th power and using \eqref{eq:Plancherel-trace-scaling}, we obtain
\begin{align*}
\|j_{q}^{\tau_{G}}(\lambda_{G}(f_{F,\varepsilon}))\|_{L^{q}(\operatorname{VN}(G),\tau_{G})}^{q}
=c_{\varepsilon}^{-1}
\tau_{G}\left(W_{\varepsilon}|\lambda_{\varepsilon}(F)|^{q}W_{\varepsilon}^{-1}\right)
&=\tau_{\varepsilon}(|\lambda_{\varepsilon}(F)|^{q})\\
&=\|j_{q}^{\tau_{\varepsilon}}(\lambda_{\varepsilon}(F))\|_{L^{q}(\operatorname{VN}(G_{\varepsilon}),\tau_{\varepsilon})}^{q}.
\end{align*}
This is exactly the asserted equality of the $\mathcal{F}L^{q}$ norms on the Fourier--Tomita core. Since the Hausdorff--Young maps $\mathcal{F}_{G,p}$ and $\mathcal{F}_{G_{\varepsilon},p}$ are continuous, density of the core in the respective $L^{p}$ spaces gives the result for every $F\in L^{p}(G_{\varepsilon})$. Notice that the input scaling $c_{\varepsilon}^{\frac{1}{p}}$, the operator scaling $c_{\varepsilon}^{-\frac{1}{q}}$, and the trace scaling in \eqref{eq:Plancherel-trace-scaling} cancel exactly.
\end{proof}

\begin{proposition}\label{prop:Cartan-contraction-lower-bound}
Under the assumptions in Theorem~\ref{thm:semisimple-sharp} and for
$1\leq p\leq2$, one has
\begin{align*}
H_{p}(G)\geq H_{p}(K\ltimes\mathfrak{p}).
\end{align*}
\end{proposition}
\begin{proof}
Let $G_{0}$ and $G_{\varepsilon}$ be given as in Lemma \ref{lem:contracted-group-law-Haar}. For $1<p<2$, $q=\frac{p}{p-1}$ and nonzero $F\in C_{c}^{\infty}(G_{0})$, it follows from the dominated convergence theorem and Lemma \ref{lem:contracted-group-law-Haar} \eqref{G0 measure convergence} that
\begin{align*}
\|F\|_{L^{p}(G_{\varepsilon},dm_{\varepsilon})}\longrightarrow\|F\|_{L^{p}(G_{0},dm_{0})},\quad\mbox{as}\quad\varepsilon\rightarrow0.
\end{align*}
Combining this with Lemma \ref{lem:exact-Plancherel-scaling}  and Lemma \ref{lem:weighted-spectral-measures},
\begin{align*}
\lim_{\varepsilon\downarrow0}
\frac{\|\mathcal{F}_{G,p}f_{F,\varepsilon}\|_{\mathcal{F}L^{q}(G)}}
{\|f_{F,\varepsilon}\|_{L^{p}(G)}}
=\lim_{\varepsilon\downarrow0}
\frac{\|\mathcal{F}_{G_{\varepsilon},p}F\|_{\mathcal{F}L^{q}(G_{\varepsilon})}}
{\|F\|_{L^{p}(G_{\varepsilon},m_{\varepsilon})}}
=\frac{\|\mathcal{F}_{G_{0},p}F\|_{\mathcal{F}L^{q}(G_{0})}}
{\|F\|_{L^{p}(G_{0},dm_{0})}},
\end{align*}
where $f_{F,\varepsilon}$ is given in Lemma \ref{lem:exact-Plancherel-scaling}.  Thus, by definition of $H_{p}(G)$ and density of $C_{c}^{\infty}(G_{0})$ in $L^{p}(G_{0})$, 
\begin{align*}
H_{p}(G)
\geq\sup_{0\neq F\in L^{p}(G_{0})}
\frac{\|\mathcal{F}_{G_{0},p}F\|_{\mathcal{F}L^{q}(G_{0})}}
{\|F\|_{L^{p}(G_{0},dm_{0})}}
=H_{p}(G_{0}).
\end{align*}

Applying Lemma \ref{lem:contracted-group-law-Haar} \eqref{G0 isomorphism}, the map $\Psi:(k,X)\longmapsto(k,\operatorname{Ad}(k)X)$ is a Haar-preserving isomorphism from $G_{0}$ to $K\ltimes_{\operatorname{Ad}}\mathfrak{p}$. If $f\in L^{p}(G_{0})$, let $\widetilde{f}=f\circ\Psi^{-1}$. The Haar-preserving group isomorphism and its induced von Neumann algebra isomorphism give that
\begin{align*}
\|\widetilde{f}\|_{L^{p}(K\ltimes\mathfrak{p})}
=\|f\|_{L^{p}(G_{0})}\quad\mbox{and}\quad
\|\mathcal{F}_{K\ltimes\mathfrak{p},p}\widetilde{f}\|_{\mathcal{F}L^{q}(K\ltimes\mathfrak{p})}
=\|\mathcal{F}_{G_{0},p}f\|_{\mathcal{F}L^{q}(G_{0})}.
\end{align*}
Therefore, 
\begin{align*}
H_{p}(G_{0})=H_{p}(K\ltimes\mathfrak{p}),
\end{align*}
and the assertion for $1<p<2$ follows.

For $p=1$ and $p=2$, Lemma~\ref{lem:H1-H2-endpoints} gives the
corresponding equalities.  The proof is complete.
\end{proof}

\subsection{The sharp constant of the Cartan motion group}\label{subsec:motion-group-constant}
The following result recovers the compact-extension result of Baklouti
and Inoue \cite{BakloutiInoue2014}; see also the classical
compact-extension theorem of Russo \cite{Russo1978}.

We use the orthogonal action, semidirect-product, and $K$-invariance
conventions from Section~\ref{subsec:preliminary-semisimple-Lie-groups}.
\begin{proposition}\label{prop:motion-group-sharp-constant}
Let $V$ be a $d$-dimensional Euclidean space on which a compact group $K$ acts orthogonally, and let the Haar measure on $K$ be normalized to have mass one. For $1\leq p\leq2$ and $q=\frac{p}{p-1}$, with $q=\infty$ when $p=1$, one has
\begin{align*}
H_{p}(K\ltimes V)=\left(p^{\frac{1}{2p}}q^{-\frac{1}{2q}}\right)^{d}.
\end{align*}
\end{proposition}
\begin{proof}
Let $1<p<2$ and $q=\frac{p}{p-1}$. For a $K$-invariant function $g\in C_{c}^{\infty}(V)$, define
\begin{align*}
(\iota g)(k,x)=g(x).
\end{align*}
Since $dk$ has mass one,
\begin{align*}
\|\iota g\|_{L^{p}(K\ltimes V)}=\|g\|_{L^{p}(V)}.
\end{align*}
For $K$-invariant $g,h$, direct calculation of the semidirect-product convolution and involution gives
\begin{align*}
\iota g*_{K\ltimes V}\iota h=\iota(g*_{V}h)\quad\mbox{and}\quad
(\iota g)^{*}=\iota(g^{*}),
\end{align*}
where $g^{*}(x)=\overline{g(-x)}$. Hence, for every integer $n\geq1$,
\begin{align*}
\tau_{K\ltimes V}\left(\left|\lambda_{K\ltimes V}(\iota g)\right|^{2n}\right)
=\left[\left((\iota g)^{*}*_{K\ltimes V}\iota g\right)^{(*_{K\ltimes V}n)}\right](e_{K},0)
&=\left[\left(g^{*}*_{V}g\right)^{(*_{V}n)}\right](0)\\
&=\tau_{V}\left(\left|\lambda_{V}(g)\right|^{2n}\right).
\end{align*}
We now justify the passage from integer moments to arbitrary real $q\geq2$.
Both $\lambda_{K\ltimes V}(\iota g)^{*}\lambda_{K\ltimes V}(\iota g)$ and $\lambda_{V}(g)^{*}\lambda_{V}(g)$ have spectrum in the common compact interval $[0,c]$, where one may take $c=\|g\|_{L^{1}(V)}^{2}$, and both have finite first trace moment. If $E_{1}$ and $E_{0}$ are their spectral measures, define finite positive measures on $[0,c]$ by
\begin{align*}
\nu_{1}(B)=\tau_{K\ltimes V}\Big(\lambda_{K\ltimes V}(\iota g)^{*}\lambda_{K\ltimes V}(\iota g)E_{1}(B)\Big)\quad\mbox{and}\quad
\nu_{0}(B)=\tau_{V}\Big(\lambda_{V}(g)^{*}\lambda_{V}(g)E_{0}(B)\Big).
\end{align*}
The preceding integer-moment identities imply
\begin{align*}
\int_{0}^{c}t^{k}\,d\nu_{1}(t)
=\tau_{K\ltimes V}\Big((\lambda_{K\ltimes V}(\iota g)^{*}\lambda_{K\ltimes V}(\iota g)\Big)^{k+1})
&=\tau_{V}\Big((\lambda_{V}(g)^{*}\lambda_{V}(g))^{k+1}\Big)\\
&=\int_{0}^{c}t^{k}\,d\nu_{0}(t),
\qquad k=0,1,2,\ldots.
\end{align*}
The density of polynomials in $C([0,c])$ gives $\nu_{1}=\nu_{0}$. Therefore,
\begin{align*}
\|\mathcal{F}_{K\ltimes V,p}(\iota g)\|_{\mathcal{F}L^{q}(K\ltimes V)}^{q}
&=\|j_{q}^{\tau_{K\ltimes V}}(\lambda_{K\ltimes V}(\iota g))\|_{L^{q}(\operatorname{VN}(K\ltimes V),\tau_{K\ltimes V})}^{q}
=\int_{0}^{c}t^{\frac{q}{2}-1}\,d\nu_{1}(t)\\
&=\int_{0}^{c}t^{\frac{q}{2}-1}\,d\nu_{0}(t)
=\|j_{q}^{\tau_{V}}(\lambda_{V}(g))\|_{L^{q}(\operatorname{VN}(V),\tau_{V})}^{q}\\
&=\|\mathcal{F}_{V,p}g\|_{\mathcal{F}L^{q}(V)}^{q}.
\end{align*}
Next, approximate a radial Euclidean Gaussian in $L^{p}(V)$ by smooth radial compactly supported functions. These functions are $K$-invariant, and Hausdorff--Young continuity permits passage to the limit. Since Gaussians attain the Euclidean Babenko--Beckner constant and the $L^{p}$ norms are preserved by $g\mapsto\iota g$, one has
\begin{align*}
H_{p}(K\ltimes V)\geq\left(p^{\frac{1}{2p}}q^{-\frac{1}{2q}}\right)^{d}.
\end{align*}
Meanwhile, since the vector subgroup $V$ is closed in $K\ltimes V$,
Theorem~\ref{thm:closed-subgroup-comparison} and
 Babenko--Beckner theorem \cite{Babenko1961,Beckner1975} on Euclidean space give that
\begin{align*}
H_{p}(K\ltimes V)\leq H_{p}(V)=\left(p^{\frac{1}{2p}}q^{-\frac{1}{2q}}\right)^{d}.
\end{align*}
Combining the last two inequalities proves the assertion for $1<p<2$.
The endpoint cases follow from Lemma~\ref{lem:H1-H2-endpoints}.
\end{proof}

Now, we give the complete proof for the semisimple Lie groups.
\begin{proof}[Proof of Theorem~\ref{thm:semisimple-sharp}]
Let $q=\frac{p}{p-1}$, with $q=\infty$ when $p=1$. Applying the Cartan decomposition and
Iwasawa decomposition, Theorem
\ref{thm:closed-subgroup-comparison}, Lemma
\ref{lem:AN-completely-solvable}, and Theorem
\ref{thm:completely-solvable-sharp},
\begin{align*}
H_{p}(G)\leq H_{p}(AN)=\left(p^{\frac{1}{2p}}q^{-\frac{1}{2q}}\right)^{\dim(AN)}
=\left(p^{\frac{1}{2p}}q^{-\frac{1}{2q}}\right)^{\dim(G/K)}.
\end{align*}
On the other hand, using Proposition \ref{prop:Cartan-contraction-lower-bound} and Proposition \ref{prop:motion-group-sharp-constant} (with $V=\mathfrak{p}$), and Lemma \ref{lem:AN-completely-solvable},
\begin{align*}
H_{p}(G)\geq H_{p}(K\ltimes\mathfrak{p})
    =\left(p^{\frac{1}{2p}}q^{-\frac{1}{2q}}\right)^{\dim\mathfrak{p}}
    =\left(p^{\frac{1}{2p}}q^{-\frac{1}{2q}}\right)^{\dim(G/K)}.
\end{align*}
The proof is complete.
\end{proof}
\begin{remark}
The proof of Theorem~\ref{thm:semisimple-sharp} includes arbitrary real rank and allows compact semisimple factors. If $G$ is compact, then $G=K$, $\dim(G/K)=0$, and the formula reduces to $H_{p}(G)=1$, as it should.
\end{remark}

\section{Examples of finite-center semisimple Lie groups}\label{sec:examples}
In this section, we give clear sharp constants of Hausdorff-Young inequality for some examples on semisimple Lie groups by using Theorem~\ref{thm:semisimple-sharp}. 
Throughout this section, $1\leq p\leq 2$ and $q=\frac{p}{p-1}$, with the usual interpretation $q=\infty$ when $p=1$.

The maximal compact subgroups and the dimension formulae used in this section are standard; see, for example, \cite{Helgason2001,Knapp2002}.
\subsection{Classical real forms}\label{subsec:classical-real-forms}

Let
\begin{align*}
 I_{r,s}=\operatorname{diag}(I_{r},-I_{s}),\quad
 J_{n}=
 \begin{pmatrix}
 0&I_{n}\\
 -I_{n}&0
 \end{pmatrix}.
\end{align*}
For complex and quaternionic matrices, $g^{*}$ denotes conjugate transpose.
The indefinite unitary and orthogonal groups are
\begin{align*}
 SU(r,s)= \left\{ g\in SL(r+s,\mathbb{C}): g^{*}I_{r,s}g=I_{r,s}\right\},\quad
 SO(r,s)=\left\{g\in SL(r+s,\mathbb{R}):g^{T}I_{r,s}g=I_{r,s}\right\}.
\end{align*}
The subscript in $SO_{0}(r,s)$ denotes the identity component of $SO(r,s)$.
Its maximal compact subgroup is locally isomorphic to
$SO(r)\times SO(s)$.  The group in the unitary table is
\begin{align*}
 S(U(r)\times U(s))
 =
 \left\{
 (u,v)\in U(r)\times U(s):
 \det(u)\det(v)=1
 \right\}.
\end{align*}

We use the convention that the real symplectic group is
\begin{align*}
 Sp(n,\mathbb{R})
 =
 \left\{
 g\in GL(2n,\mathbb{R}):
 g^{T}J_{n}g=J_{n}
 \right\}.
\end{align*}
The compact group $Sp(n)$ is the quaternionic unitary group
\begin{align*}
 Sp(n)
 =
 \left\{
 g\in GL(n,\mathbb{H}):
 g^{*}g=I_{n}
 \right\}.
\end{align*}
For $r,s\geq1$, the group $Sp(r,s)$ consists of quaternionic linear
transformations preserving the quaternionic Hermitian form of signature
$(r,s)$:
\begin{align*}
 Sp(r,s)
 =
 \left\{
 g\in GL(r+s,\mathbb{H}):
 g^{*}I_{r,s}g=I_{r,s}
 \right\}.
\end{align*}
Its maximal compact subgroup is $Sp(r)\times Sp(s)$.

The symbols $SU^{*}(2n)$ and $SO^{*}(2n)$ denote quaternionic real forms.
If the quaternionic structure on $\mathbb{C}^{2n}$ is
$v\longmapsto J_{n}\overline{v}$, then they may be realized as
\begin{align*}
 SU^{*}(2n)
 &=
 \left\{
 g\in SL(2n,\mathbb{C}):
 gJ_{n}=J_{n}\overline{g}
 \right\}
 \simeq SL(n,\mathbb{H}),\\
 SO^{*}(2n)
 &=
 \left\{
 g\in SO(2n,\mathbb{C}):
 gJ_{n}=J_{n}\overline{g}
 \right\}.
\end{align*}
Their maximal compact subgroups are $Sp(n)$ and $U(n)$, respectively.  The
star in $SO^{*}(2n)$ and $SU^{*}(2n)$ is part of Cartan's notation for the
real form; it does not denote the adjoint of a matrix.

The basic compact-group dimensions are
\begin{align*}
 \dim SO(n)&=\frac{n(n-1)}{2},
 &\dim U(n)&=n^{2},\\
 \dim SU(n)&=n^{2}-1,
 &\dim Sp(n)&=n(2n+1).
\end{align*}
Subtracting $\dim K$ from $\dim G$ gives the following values.

\begin{table}[htbp]
\centering
\small
\begin{tabular}{@{}>{\raggedright\arraybackslash}p{0.25\textwidth}p{0.22\textwidth}p{0.18\textwidth}p{0.23\textwidth}@{}}
\hline
$G$ & maximal compact subgroup $K$ & $\dim(G/K)$ & $H_{p}(G)$ \\
\hline
$SL(n,\mathbb{R})$, $n\geq2$
 & $SO(n)$
 & $\frac{(n-1)(n+2)}{2}$
 & $B_{p}^{\frac{(n-1)(n+2)}{2}}$ \\
$SL(n,\mathbb{C})$, as a real group
 & $SU(n)$
 & $n^{2}-1$
 & $B_{p}^{n^{2}-1}$ \\
$SU(r,s)$, $r,s\geq1$
 & $S(U(r)\times U(s))$
 & $2rs$
 & $B_{p}^{2rs}$ \\
$SO_{0}(r,s)$, $r,s\geq1$, \mbox{$r+s\geq3$}
 & $SO(r)\times SO(s)$
 & $rs$
 & $B_{p}^{rs}$ \\
\hline
\end{tabular}
\caption{Special linear, unitary, and orthogonal families.}
\label{tab:classical-examples-one}
\end{table}

\begin{table}[htbp]
\centering
\small
\begin{tabular}{@{}>{\raggedright\arraybackslash}p{0.25\textwidth}p{0.22\textwidth}p{0.18\textwidth}p{0.23\textwidth}@{}}
\hline
$G$ & maximal compact subgroup $K$ & $\dim(G/K)$ & $H_{p}(G)$ \\
\hline
$Sp(n,\mathbb{R})$, $n\geq1$
 & $U(n)$
 & $n(n+1)$
 & $B_{p}^{n(n+1)}$ \\
$Sp(r,s)$, $r,s\geq1$
 & $Sp(r)\times Sp(s)$
 & $4rs$
 & $B_{p}^{4rs}$ \\
$SO^{*}(2n)$, $n\geq2$
 & $U(n)$
 & $n(n-1)$
 & $B_{p}^{n(n-1)}$ \\
$SU^{*}(2n)\simeq SL(n,\mathbb{H})$, \mbox{$n\geq2$}
 & $Sp(n)$
 & $(n-1)(2n+1)$
 & $B_{p}^{(n-1)(2n+1)}$ \\
\hline
\end{tabular}
\caption{Symplectic and quaternionic families and the remaining classical real forms.}
\label{tab:classical-examples-two}
\end{table}
The restrictions on the parameters in Tables~\ref{tab:classical-examples-one} and \ref{tab:classical-examples-two} exclude the non-semisimple low-dimensional cases and restrict the displayed families to noncompact real forms.  In particular, $SO^{*}(4)$ is allowed: it is semisimple but not simple.
The omitted group $SU^{*}(2)\simeq SU(2)$ is compact and therefore has Hausdorff--Young constant $1$.

The group $SO_{0}(r,s)$ has finite center.  Its connected spin double cover $\operatorname{Spin}_{0}(r,s)\longrightarrow SO_{0}(r,s)$ has finite kernel and finite center, and its maximal compact subgroup is the
inverse image of $SO(r)\times SO(s)$, locally isomorphic to $\operatorname{Spin}(r)\times\operatorname{Spin}(s)$.  Hence this spin group---and, by the preceding global-form observation, every connected finite-center group locally isomorphic to it---satisfies
\begin{align*}
 H_{p}\bigl(\operatorname{Spin}_{0}(r,s)\bigr)=B_{p}^{rs}.
\end{align*}
This statement deliberately does not include an infinite-sheeted universal cover with infinite center.

\subsection{Complex and compact groups}\label{subsec:complex-compact-examples}

The symbols
\begin{align*}
 G_{2}(\mathbb{C}),\quad
 F_{4}(\mathbb{C}),\quad
 E_{6}(\mathbb{C}),\quad
 E_{7}(\mathbb{C}),\quad
 E_{8}(\mathbb{C})
\end{align*}
denote connected complex simple groups whose complex Lie algebras have the indicated Dynkin types.  A chosen global form is always assumed to have finite center.  The symbols $G_{2},$ $F_{4},$ $E_{6},$ $E_{7},$ $E_{8}$
without an argument denote compact real forms when they occur as compact groups below.

Let $G_{\mathbb{C}}$ be a connected complex semisimple group with finite center, regarded as a real Lie group, and let $K$ be a compact real form. Then
\begin{align*}
 \dim_{\mathbb{R}}G_{\mathbb{C}}=2\dim_{\mathbb{R}}K,
 \qquad\dim_{\mathbb{R}}(G_{\mathbb{C}}/K)=\dim_{\mathbb{C}}G_{\mathbb{C}}.
\end{align*}
By Theorem \ref{thm:semisimple-sharp}, it follows that
\begin{align}
 H_{p}(G_{\mathbb{C}})=B_{p}^{\dim_{\mathbb{C}}G_{\mathbb{C}}}.
 \label{eq:complex-semisimple-examples}
\end{align}
In particular,
\begin{align*}
 H_{p}(SO(n,\mathbb{C}))=B_{p}^{\frac{n(n-1)}{2}},\quad  n\geq3,\quad\mbox{and}\quad
 H_{p}(Sp(n,\mathbb{C}))=B_{p}^{n(2n+1)},\quad  n\geq1,
\end{align*}
where
\begin{align*}
 Sp(n,\mathbb{C})
 =
 \left\{
 g\in GL(2n,\mathbb{C}):
 g^{T}J_{n}g=J_{n}
 \right\}.
\end{align*}
The groups are regarded as real Lie groups and their maximal compact subgroups
are $SO(n)$ and $Sp(n)$, respectively.  Formula
\eqref{eq:complex-semisimple-examples} gives
\begin{align*}
 \begin{array}{c|ccccc}
 G_{\mathbb{C}} & G_{2}(\mathbb{C}) & F_{4}(\mathbb{C}) & E_{6}(\mathbb{C})
 & E_{7}(\mathbb{C}) & E_{8}(\mathbb{C})\\
 \hline
 H_{p}(G_{\mathbb{C}}) & B_{p}^{14} & B_{p}^{52} & B_{p}^{78}
 & B_{p}^{133} & B_{p}^{248}.
 \end{array}
\end{align*}

At the opposite extreme, if $G$ is compact, then one may take $K=G$. Thus $H_{p}(G)=1.$ For example, $H_{p}(Sp(n))=1$.  More generally, compact semisimple factors do not contribute to the exponent in Theorem~\ref{thm:semisimple-sharp}.

\subsection{Selected exceptional real forms}\label{subsec:exceptional-real-form-examples}

Let $\mathfrak{g}$ be a real form of a complex exceptional simple Lie algebra
and write its Cartan decomposition as $\mathfrak{g}=\mathfrak{k}\oplus\mathfrak{p}.$
The Killing form is negative definite on $\mathfrak{k}$ and positive definite
on $\mathfrak{p}$.  The Cartan signature, also called the character of the
real form, is
\begin{align*}
 s(\mathfrak{g})
 =
 \dim(\mathfrak{p})-\dim(\mathfrak{k}).
\end{align*}
The notation $G_{2(2)}$, $F_{4(-20)}$, $E_{6(-26)}$, and so forth denotes a
connected finite-center group whose Lie algebra is the real form of the
indicated complex Dynkin type with Cartan signature equal to the integer in
parentheses.  Thus, in $E_{6(-26)}$, the index $6$ belongs to the Dynkin symbol
$E_{6}$, whereas $-26$ identifies the real form; the latter is neither its
dimension nor its real rank.

The groups
\begin{align*}
 G_{2(2)},\quad
 F_{4(4)},\quad
 E_{6(6)},\quad
 E_{7(7)},\quad
 E_{8(8)}
\end{align*}
are the split real forms.  For example, $G_{2(2)}$ is the automorphism group
of the split octonion algebra.  The groups
\begin{align*}
 F_{4(-20)},\quad
 E_{6(-26)},\quad
 E_{7(-25)},\quad
 E_{8(-24)}
\end{align*}
are selected nonsplit noncompact real forms; in particular,
$F_{4(-20)}/\operatorname{Spin}(9)$ is the Cayley hyperbolic plane.  In the
maximal compact subgroup column, symbols such as $F_{4}$, $E_{6}$, and $E_{7}$
denote compact real forms.  Products and spin groups in that column are
understood up to the finite central quotient required by the chosen global
form.

For completeness, Table~\ref{tab:exceptional-real-forms} records their
symmetric-space dimensions and sharp constants.

\begin{table}[htbp]
\centering
\small
\begin{tabular}{@{}p{0.25\textwidth}p{0.22\textwidth}p{0.18\textwidth}p{0.23\textwidth}@{}}
\hline
$G$ & maximal compact subgroup $K$ & $\dim(G/K)$ & $H_{p}(G)$ \\
\hline
$G_{2(2)}$ & $SO(4)$ & $8$ & $B_{p}^{8}$ \\
$F_{4(4)}$ & $Sp(3)\times Sp(1)$ & $28$ & $B_{p}^{28}$ \\
$F_{4(-20)}$ & $\operatorname{Spin}(9)$ & $16$ & $B_{p}^{16}$ \\
$E_{6(6)}$ & $Sp(4)$ & $42$ & $B_{p}^{42}$ \\
$E_{6(-26)}$ & $F_{4}$ & $26$ & $B_{p}^{26}$ \\
$E_{7(7)}$ & $SU(8)$ & $70$ & $B_{p}^{70}$ \\
$E_{7(-25)}$ & $E_{6}\times U(1)$ & $54$ & $B_{p}^{54}$ \\
$E_{8(8)}$ & $\operatorname{Spin}(16)$ & $128$ & $B_{p}^{128}$ \\
$E_{8(-24)}$ & $E_{7}\times SU(2)$ & $112$ & $B_{p}^{112}$ \\
\hline
\end{tabular}
\caption{Selected exceptional real forms.}
\label{tab:exceptional-real-forms}
\end{table}

For example,
\begin{align*}
 \dim E_{7(-25)}-\dim(E_{6}\times U(1))
 =133-(78+1)=54,
\end{align*}
and
\begin{align*}
 \dim E_{8(-24)}-\dim(E_{7}\times SU(2))
 =248-(133+3)=112.
\end{align*}
The compact exceptional real forms again have quotient dimension $0$ and hence sharp constant $1$.

\end{document}